\documentclass[english,reqno]{amsart}
\usepackage[T1]{fontenc}
\usepackage[utf8]{inputenc}
\usepackage{color}
\usepackage{babel}
\usepackage{amstext}
\usepackage{amsthm}
\usepackage{amssymb}
\usepackage{geometry}
\usepackage[pdfusetitle,
 bookmarks=true,bookmarksnumbered=false,bookmarksopen=false,
 breaklinks=false,pdfborder={0 0 0},pdfborderstyle={},backref=false,colorlinks=true]
 {hyperref}
\hypersetup{
 linkcolor=brown, citecolor=blue, pdfstartview={FitH}, hyperfootnotes=false, unicode=true}

\makeatletter
\numberwithin{figure}{section}
\theoremstyle{plain}
\newtheorem{thm}{\protect\theoremname}[section]
\theoremstyle{definition}
\newtheorem{defn}[thm]{\protect\definitionname}
\theoremstyle{plain}
\newtheorem{prop}[thm]{\protect\propositionname}
\theoremstyle{plain}
\newtheorem{lem}[thm]{\protect\lemmaname}
\theoremstyle{definition}
\newtheorem{example}[thm]{\protect\examplename}
\theoremstyle{plain}
\newtheorem{cor}[thm]{\protect\corollaryname}
\theoremstyle{plain}
\newtheorem{fact}[thm]{\protect\factname}
\theoremstyle{definition}
\newtheorem{problem}[thm]{\protect\problemname}

\usepackage{tikz-cd}
\usepackage{enumitem}
\usepackage{bbm}
\usepackage{centernot}
\usepackage{wasysym}
\setlist[enumerate]{itemsep=3pt, topsep=5pt}
\setlist[itemize]{itemsep=3pt, topsep=5pt}

\newtheorem{theoremalpha}{Theorem}

\makeatother

\providecommand{\corollaryname}{Corollary}
\providecommand{\definitionname}{Definition}
\providecommand{\examplename}{Example}
\providecommand{\factname}{Fact}
\providecommand{\lemmaname}{Lemma}
\providecommand{\problemname}{Problem}
\providecommand{\propositionname}{Proposition}
\providecommand{\theoremname}{Theorem}

\begin{document}
\title[Effective bases and effective second countability]{Effective bases and notions of effective second countability in computable
analysis }
\author{Vasco Brattka}
\address{Fakultät für Informatik, Universität der Bundeswehr München, Werner-Heisenberg-Weg
39, 85577 Neubiberg, Germany }
\address{Department of Mathematics and Applied Mathematics, University of Cape
Town, Private Bag X3, Rondebosch 7701, South Africa}
\email{vasco.brattka@cca-net.de}
\author{Emmanuel Rauzy}
\address{LACL, Université Paris-Est Créteil, 61 avenue du Général de Gaulle,
94010 Créteil}
\email{emmanuel.rauzy@u-pec.fr }
\begin{abstract}
We compare the approaches of Schröder and Weihrauch to computable
topology. These approaches are both based on the theory of representations,
but they focus on a priori different representations of open sets.
We show that the two approaches are compatible. While it is well known
that Schröder’s approach is more general because it can be applied
to non second countable spaces, we insist on the fact that it is also
more general even when dealing only with second countable spaces,
because a represented space can be second countable without being
\emph{effectively second countable}. 

We revisit Schröder's Effective Metrization Theorem, by showing that
it characterizes those represented spaces that embed into computable
metric spaces: those are the computably second countable strongly
computably regular represented spaces. 

Finally, we study different forms of open choice problems. We show
that having a computable open choice is equivalent to being computably
separable, but that the ``non-total open choice problem'', i.e.,
open choice restricted to open sets that have non-empty complement,
interacts with effective second countability in a satisfying way. 
\end{abstract}

\maketitle
\tableofcontents{}

\section{Introduction}

Choosing a good framework and correct definitions is one of the most
important and sometimes most challenging problems in computable mathematics.
Indeed, it is often the case that once a good framework has been found,
many results that seemed non-trivial before can easily be proven.
To see that this is the case in computable analysis, one can for instance
look at early articles of Markov \cite{Markov1963} and Lacombe \cite{Lacombe1957b},
which translated to today's vocabulary end up having very little content,
and compare them to Pauly's 2016 survey paper \cite{Pauly2016} which
deals with a great many notions and results in a mere 22 pages. 

$ $

We are concerned here with computable topology as studied in Weihrauch's
German school of computable analysis, we follow a ``Type 2 approach'':
the basic notion of computable function is set on Baire space, and
it is transferred from Baire space to other spaces with cardinality
that of the continuum thanks to partial surjections which are called
\emph{representations}. 

Our starting point is that within Weihrauch's school of computable
analysis two approaches to studying computable topology coexist. 

Namely: 
\begin{itemize}
\item The approach devised by Matthias Schröder in his PhD thesis \cite{Schroeder2002}
(see also \cite{Schroeder2001,Schroeder2021}), which uses the representation
of the Sierpiński space to define a representation of open sets on
any represented space equipped with the final topology of its representation.
\item The approach that Weihrauch developed with Grubba in \cite{Weihrauch2009ElementaryCT},
which relies on a numbered base $(B_{i})_{i\in\mathbb{N}}$ for a
space $X$ to define a representations of $X$ and a representation
of the open sets of $X$.
\end{itemize}
Ideas very similar to those of Schröder were developed independently
by Paul Taylor \cite{Taylor2011} and Martín Escardó \cite{Escardo2004},
outside of the framework of represented sets. 

It is clear that Schröder's approach is more general than Weihrauch's,
in that it applies to spaces that need not be second countable. Yet
the exact relationship between these two approaches had not been clarified
until now. 

In the present paper, we explain the following: 
\begin{itemize}
\item The approaches of Schröder and Weihrauch to computable topologies
are compatible. 
\item The approach of Schröder is strictly more general than Weihrauch's,
\emph{even for second countable spaces}: not only does it encompass
more topological spaces, but even on the topological spaces that fall
under both approaches, purely computability theoretical phenomena
occur in Schröder's approach that are not accounted for in Weihrauch's
approach. 
\item Weihrauch and Grubba's notion of ``computable topological space''
in fact corresponds to a form of \emph{effective second countability.}
\item This notion of effective second countability is often the most useful
one, yet other weaker notions can be relevant, because the ``effectivization''
of some theorems can already be obtained thanks to those weaker notions. 
\end{itemize}
$ $

Another way to understand the present article is as follows. 

One of the key notions to study topology thanks to the theory of represented
spaces is that of \emph{admissible representation}, because continuous
maps between admissibly represented spaces always have a continuous
realizer. Also, when studying computable topology, \emph{computably
admissible representations} play a central role, because these are
the representations that guarantee the equivalence between ``being
computable'' and ``being effectively continuous''. See Section
\ref{sec:Preliminaries} for more details. 

One of the earliest representations ever studied is the \emph{standard
representation} introduced by Kreitz and Weihrauch in \cite{Kreitz1985}.
Given a second countable $\text{T}_{0}$ space $X$ equipped with
a countable subbase $(B_{i})_{i\in\mathbb{N}}$, this representation
is given by 
\[
\rho(f)=x\iff\text{Im}(f)-1=\{n\in\mathbb{N},\,x\in B_{n}\},
\]
where $\text{Im}(f)-1=\{n\in\mathbb{N},\,\exists k\in\mathbb{N},n+1=f(k)\}$.
Note that every totally numbered subbase of $X$ induces a standard
representation, and thus each second countable space admits multiple
standard representations. However, Kreitz and Weihrauch have shown
that, up to continuous equivalence of representations, the choice
of a standard representation is inconsequential:
\begin{thm}
[\cite{Kreitz1985}]All standard representations of a second countable
space are continuously equivalent. Every admissible representation
of a second countable space is continuously equivalent to a standard
representation. 
\end{thm}

It is clear that not all standard representations of a second countable
space need to be computably equivalent. What we want to emphasize
here is that the effective analogue of the second point in the above
theorem also fails: \emph{a computably admissible representation of
a second countable space does not have to be computably equivalent
to a standard representation.} 

And thus, while fixing a standard representation on a second countable
space can be seen as a ``neutral operation'' from the point of view
of topology, this fact does not carry over to the study of computability:
in a context where representations are considered up to computable
translations, studying only spaces that are equipped with standard
representations amounts, implicitly, to imposing certain computability-theoretic
assumptions. The representations that are computably equivalent to
a standard representation are exactly those that give rise to what
we call computably second countable spaces. 

In this paper, we show that the notion of a computably second countable
represented space is extremely robust, as it emerges from a wide variety
of different approaches. 

But we also show that among representations that are not computably
equivalent to a standard representation, there is a range of weaker
notions of effective second countability that can be, in different
contexts, relevant. 

\subsection{Notions of effective bases and their relations }

In Section \ref{sec:Notions-of-bases}, we introduce in detail the
notions of bases that we consider. Here we only sketch the definitions.
\begin{itemize}
\item \textbf{Semi-effective bases.} A semi-effective base is a set of uniformly
open sets that form a base. These sets are ``constructively open'',
but the assumption that they form a base is purely classical. 
\item \textbf{Lacombe bases.} This is the notion of base that follows from
the classical statement ``a set $\mathfrak{B}$ forms a base for
a topological space $X$ if the open sets are exactly the sets that
can be written as unions of elements of $\mathfrak{B}$''. The effective
version of this statement will say that ``open sets can uniformly
be written as overt unions of basic sets''. This notion is named
after Lacombe following the article \cite{Lacombe1957}. Papers that
have used this approach include \cite{Hoyrup2016,Amir2023,Korovina2008,KOROVINA2016,Grubba2007,Wei08,Hoyrup2019},
for the second countable case, and \cite{Bauer2012} for a more general
setting. 
\item \textbf{Nogina bases.} This is the notion of base that follows from
the classical statement ``a set $\mathfrak{B}$ forms a base for
a topological space $X$ if a set $O$ is open if and only if for
any $x\in O$ there is $B\in\mathfrak{B}$ with $x\in B\subseteq O$''.
Nogina was the first to define a notion a base by using an effective
version of this statement \cite{Nogina1966,Nogina_1969}. Recent use
of it can be found in \cite{GREGORIADES2016}. 
\item \textbf{Representation subbases.} A representation subbase of a set
$X$ is a subbase which is used to define a representation of $X$,
by saying that the name of a point should encode the characteristic
function of the set of basic sets to which it belongs. Here, characteristic
functions are considered to have the Sierpiński space as their codomain,
rather than the discrete two-point space $\{0,1\}$. This is one of
two possible generalizations of the standard representation considered
by Weihrauch and Kreitz \cite{Kreitz1985}. It appears for instance
in \cite{Bauer2025}. 
\item \textbf{Enumeration subbase.} This is a second definition where a
subbase is used to define a representation. Here the name of a point
is a countable list of names of basic sets that define a \emph{formal
neighborhood base }of this point. This notion of base applies only
to second countable spaces, but possibly to uncountable bases of second
countable spaces. This definition originates in \cite{Weihrauch1987},
where Weihrauch defined representations by describing points thanks
to neighborhood bases. Following Spreen \cite{Spreen2001}, it was
explained in \cite{Rauzy2025} how relying on a notion of formal inclusion
relation and the induced notion of formal neighborhood base yields
a more robust definition. 
\end{itemize}
To each notion of computable base, we associate a notion of computable
second countability, following the following definition scheme:
\begin{defn}
A represented space $(X,\rho)$ is \emph{{[}Nogina, Lacombe, etc{]}
second countable} if it admits a totally numbered {[}Nogina, Lacombe,
etc{]} base, i.e., if there exists a sequence $(B_{i})_{i\in\mathbb{N}}$
which forms a {[}Nogina, Lacombe, etc{]} base of $(X,\rho)$.
\end{defn}

Note that we allow for some non-admissible and non-computably admissible
representations. Indeed, some notions of computable bases automatically
imply effective admissibility of the representation under scrutiny,
while others do not, we can distinguish between them only because
we are not restricting our attention to computably admissible representations
in the first place. We abbreviate ``computably admissible'' by $\text{CT}_{0}$. 

\begin{figure}
 	\begin{center} 		
\begin{tikzcd}[row sep=1.2cm, column sep=1.1cm]  			& \text{Classical base} &\\  		
	& \text{Semi-effective base} \arrow[u] &\\ 
	 			\text{Nogina base} \arrow[ur] &\text{Lacombe base }\arrow[u] 	\arrow[l,  "\tiny{\text{+local overt choice}}"]& \text{Representation subbase} \arrow[ul] \\  		
	\text{Enumeration subbase} \arrow[u] &\text{Lacombe and CT$_0$} \arrow[u] \arrow[ru, bend right=5]  	
			& \\   	
	\end{tikzcd} 		\vspace{-8mm} 	\end{center} 

\caption{\label{fig:Notions-of-bases}Notions of bases}
\end{figure}
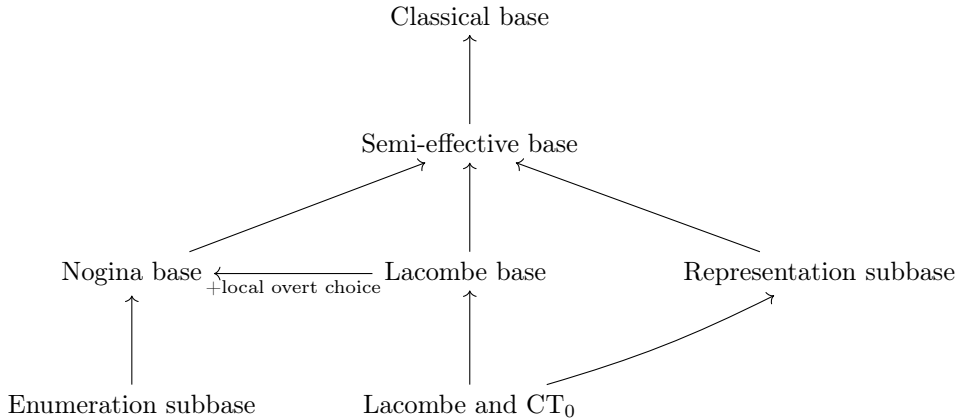

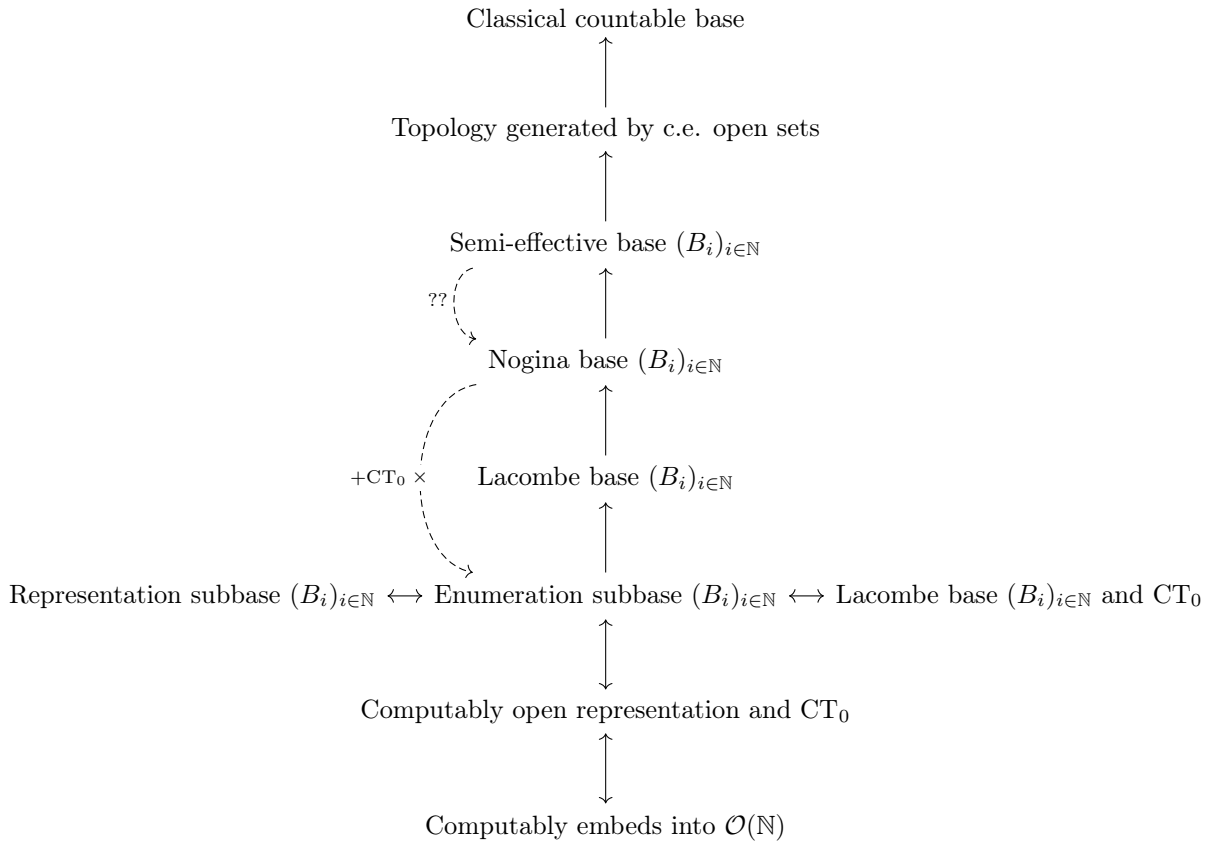
\begin{figure}
\bigskip \bigskip \bigskip
\begin{center}
\begin{tikzcd}[row sep=large, column sep=-0.5cm] 	& \text{Classical countable base} &\\ 
	& \text{Topology generated by c.e. open sets} \arrow[u] &\\ 
	& \text{Semi-effective base $(B_i)_{i\in\mathbb{N}}$} \arrow[u]\arrow[d, bend right=80, "??" swap, pos=0.45, dashed]  &\\ 
	&\text{Nogina base $(B_i)_{i\in\mathbb{N}}$} \arrow[u] \arrow[dd, bend right=80, "+ \text{CT}_0\,\,"', "\times" description, dashed] & \\ 	&	\text{Lacombe base $(B_i)_{i\in\mathbb{N}}$} \arrow[u] & \\ 
	\text{Representation subbase $(B_i)_{i\in\mathbb{N}}$} \arrow[r, leftrightarrow] &  	\text{Enumeration subbase $(B_i)_{i\in\mathbb{N}}$} \arrow[u] & 	\text{Lacombe base $(B_i)_{i\in\mathbb{N}}$ and CT$_0$} \arrow[l, leftrightarrow]  \\ 
	 &\text{Computably open representation and CT$_0$}\arrow[u, leftrightarrow] & \\ &\text{Computably embeds into $\mathcal{O}(\mathbb{N})$} \arrow[u, leftrightarrow] & \\ \end{tikzcd}
\vspace{-5mm}
\end{center}\caption{\label{fig:Notions-of-effective}Notions of computable second countability}
\end{figure}

Our main theorem is summarized in Figures \ref{fig:Notions-of-bases}
and \ref{fig:Notions-of-effective}. Because Figure \ref{fig:Notions-of-bases}
expresses, in a concise way, 8 implications, and Figure \ref{fig:Notions-of-effective}
13 implications and a ``conjectured non-implication'', it seems
reasonable to leave the figures as a statement of our theorem, instead
of listing all these implications explicitly. 

\begin{theoremalpha}	\label{thm:A-All-the-implications between bases}All
the implications between the different notions of computable bases
and of computable second countability appear in Figure \ref{fig:Notions-of-bases}
and \ref{fig:Notions-of-effective}. 

\end{theoremalpha}

Note that the two figures express something different: when comparing
notions of bases (Figure \ref{fig:Notions-of-bases}), we ask ``is
it the case that any base of this type is also automatically a base
in that sense''. (In the case of subbases we allow ourselves to replace
a subbase by the base it induces.) On the other hand, when comparing
notions of second countability (Figure \ref{fig:Notions-of-effective}),
the question we consider is: ``must a space that admits a base in
this sense also admit a base in that sense''. Those bases can, a
priori, be very different from one another. 

The five equivalent notions that appear at the bottom of Figure \ref{fig:Notions-of-effective}
define what we call\emph{ computable second countability}. Note that
three of the characterizations we present were obtained independently
in \cite{Neu25}. 

If $(X,\rho)$ is a represented space, and $Y$ is a subset of $X$,
there are two natural representations of the open sets of $Y$: one
can first restrict $\rho$ to $Y$, and then take the associated Sierpiński
representation, or first take the Sierpiński representation for $X$,
and consider the trace of this representation on $Y$. If these two
representations agree, $Y$ is called \emph{a computably sequential}
subset of $X$ (this notion was introduced by Bauer in the context
of synthetic topology under the name \emph{intrinsic subset} \cite{Bauer2025}).
An important feature of computable second countability is the following: 

\begin{theoremalpha}	\label{thm:Intrinsic Embedding-1}Let $(X,\rho)$
be a computably second countable represented space, let $Y\subseteq X$
be a subset of $X$, and equip it with the induced representation
$\rho_{\vert Y}$. Then $(Y,\rho_{\vert Y})$ is also computably second
countable, and it is a computably sequential subset of $(X,\rho)$.

\end{theoremalpha}

\subsection{More on the Grubba-Weihrauch approach }

Let us quote the Grubba-Weihrauch definition of a computable topological
space that appears in \cite{Weihrauch2009ElementaryCT}. Denote by
$W_{i}=\text{dom}(\varphi_{i})$ the usual numbering of c.e.\ subsets
of $\mathbb{N}$. 
\begin{defn}
[\cite{Weihrauch2009ElementaryCT}]\label{def:Usual Definition }A
\emph{computable topological space} is a pair $(X,(B_{i})_{i\in\mathbb{N}})$,
where $X$ is a set and $(B_{i})_{i\in\mathbb{N}}$ is the base of
a $\text{T}_{0}$ topology on $X$ for which there exists a computable
function $f:\mathbb{N}^{2}\rightarrow\mathbb{N}$ such that for any
$i$, $j$ in $\mathbb{N}$: 
\[
B_{i}\cap B_{j}=\underset{k\in W_{f(i,j)}}{\bigcup}B_{k}.
\]
\end{defn}

In fact, the above definition is systematically studied together with
two representations that are induced by the base $(B_{i})_{i\in\mathbb{N}}$:
a representation of points and a representation of open sets. 
\begin{defn}
[\cite{Weihrauch2009ElementaryCT}]\label{def:representations induced by Lacombe basis }Let
$(X,(B_{i})_{i\in\mathbb{N}})$ be a computable topological space
as above. Define a representation $\theta^{+}:\subseteq\mathbb{N}^{\mathbb{N}}\rightarrow\mathcal{O}(X)$
of the open sets of $X$ by: 
\[
\theta^{+}(f)=\underset{\{n,\,\exists p\in\mathbb{N},\,f(p)=n+1\}}{\bigcup}B_{n}.
\]

Define also a representation $\rho$ of $X$ by the following formula:
\[
\rho(f)=x\iff\text{Im}(f)=\{n\in\mathbb{N},\,x\in B_{n}\}.
\]
\end{defn}

\begin{figure}
	\begin{center} 		\begin{tikzcd}[row sep=large, column sep=3cm] 			& |[alias=O]|  \text{Representation of open sets} \\  			\text{Base}   			\arrow[ur, bend left=6, "\text{Lacombe base}" sloped,  to=O.west]  			\arrow[dr, bend right=6, "\text{Representation base}" sloped, to=1.west] &  \\  			& |[alias=1]|  \text{Representation of points}   \\  		\end{tikzcd}
\end{center}\caption{\label{fig:Weihrauch-Grubba-approach}Weihrauch-Grubba approach}
\end{figure}
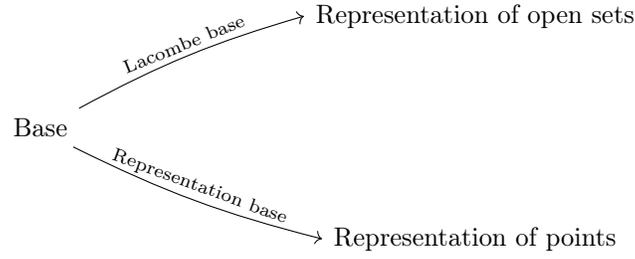

\begin{figure}
 	\begin{center}	\begin{tikzcd}[row sep=large, column sep=3cm] 		& |[alias=O]|  \text{ Representation of open sets} \\ 		\text{Base}  \arrow[dr,bend right=6,  "\text{Representation base}" sloped,  to=1.west]&  \\ 		&|[alias=1]|  \text{Representation of points }\arrow[uu,bend right=15, "\text{Sierpi\'nski representation}", swap]  \\ 	\end{tikzcd}\end{center}\caption{\label{fig:Sierpi=000144ski-representation-in schroeder}Sierpiński
representation in Schröder's approach}
\end{figure}
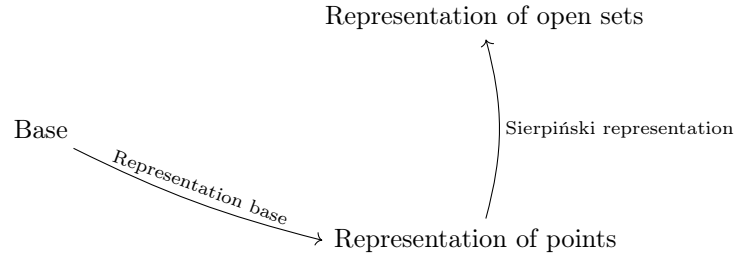

Thus in the Weihrauch-Grubba approach, a numbered basis is used to
define both a representation of open sets and a representations of
points. This is not the case in Schröder's approach to computable
topology, where the representation of open sets that is considered
is always the one associated to the Sierpiński representation. A priori,
this could lead to some conflicts, if the two approaches focus on
non-equivalent representations of open sets, but we show that this
can never arise. We summarize this in Figures \ref{fig:Weihrauch-Grubba-approach},
\ref{fig:Sierpi=000144ski-representation-in schroeder} and \ref{fig:Compatibility-of-the approaches}. 
\begin{itemize}
\item In the Weihrauch-Grubba approach, a base is used to define two representations,
as represented on Figure \ref{fig:Weihrauch-Grubba-approach}. 
\end{itemize}

\begin{itemize}
\item Following Schröder's approach, the representation of points automatically
induces a representation of open sets: this is represented on Figure
\ref{fig:Sierpi=000144ski-representation-in schroeder}. 
\end{itemize}

\begin{itemize}
\item We then express the fact that the Weihrauch-Grubba and Schröder approaches
are compatible via a commutative diagram represented on Figure \ref{fig:Compatibility-of-the approaches}. 
\end{itemize}
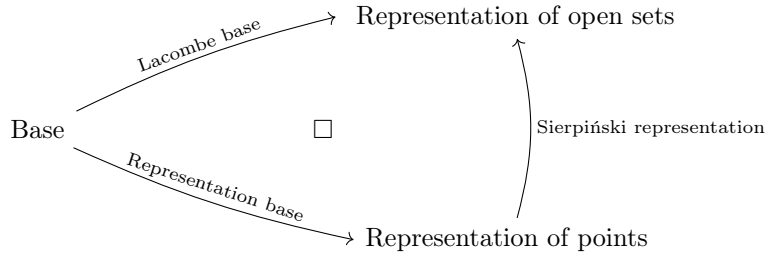
\begin{figure}
	\begin{center} 	\begin{tikzcd}[row sep=large, column sep=3cm] 		& |[alias=O]|  \text{ Representation of open sets} \\ 		\text{Base}\arrow[ur, bend left=6, "\text{Lacombe base}" sloped,  to=O.west]    \arrow[dr,bend right=6,  "\text{Representation base}" sloped,  to=1.west]&  \square \hspace{5cm} \\ 		&|[alias=1]|  \text{Representation of points }\arrow[uu,bend right=15, "\text{Sierpi\'nski representation}", swap]  \\ 	\end{tikzcd} 	\end{center}\caption{\label{fig:Compatibility-of-the approaches}Compatibility of the approaches}
\end{figure}

\medskip

Matthew de Brecht has investigated the Grubba-Weihrauch notion of
computable topological space in \cite{Brecht2020a}. He proved in
particular that each computable function $f$ that computes intersections
for a certain base, i.e., the function $f$ that appears in 
\[
B_{i}\cap B_{j}=\underset{k\in W_{f(i,j)}}{\bigcup}B_{k},
\]
can be associated to a unique maximal topological space which is a
precomputable quasi-Polish space \cite{Brecht2020}. This provides
a method for uniquely specifying a topological space and an associated
representation thanks to a finite amount of data. This approach is
orthogonal to the one presented here: we show that Schröder's approach
leads to the same representations as those introduced thanks to Definition
\ref{def:representations induced by Lacombe basis }, but foregoing
the function $f$ entirely, de Brecht has shown what is exactly the
information that can be recovered from $f$. 

\subsection{Refining results that rely on effective second countability}

Consider a classical mathematical theorem of the form:
\begin{equation}
\text{If a topological set \ensuremath{X} is second countable and A, then B. }\label{eq:1}
\end{equation}

Thanks to the different notions of effective second countability that
we consider, the following problem now makes sense: which one of the
effective second countability notions is the weakest sufficient notion
to establish an effective version of this classical result?

Many results obtained in the Weihrauch-Grubba formalism could be revisited
with these questions in mind, in particular the results of \cite{Weihrauch2010,Weihrauch2013}. 

In the present article, we are content with providing two examples
of classical statements of the form \eqref{eq:1} that require a strong
and a weak version of ``effective second countability'' respectively.
These two examples seem to us particularly important. 

The first example is Matthias Schröder's Effective Metrization Theorem
\cite{Schroeder1998}. Using ideas of Amir and Hoyrup \cite{Amir2023},
we prove that this theorem is sharp, in that it can be written with
an ``if and only if'' statement. 

\begin{theoremalpha}	[Schröder-Urysohn Effective Metrization]The
following are equivalent for a represented space $(X,\rho)$:
\begin{enumerate}
\item $(X,\rho)$ computably embeds into the Hilbert cube,
\item $(X,\rho)$ is computably second countable and strongly computably
regular.
\end{enumerate}
The above imply, but are not equivalent to: 
\begin{enumerate}
\item [(3)]\setcounter{enumi}{3} $(X,\rho)$ is computably second countable
and admits a computable metric that generates its topology. 
\end{enumerate}
\end{theoremalpha}

Previous versions of the above theorem stated (2)$\implies$(3). The
most important implication is actually (2)$\implies$(1). This implication
cannot be recovered from (2)$\implies$(3) since (3) is a strictly
weaker statement. 

We note that, modulo Schröder's results from \cite{Schroeder1998},
the proof of the above theorem is very easy. The reason why we emphasize
it so is that it seemed to us very important to understand precisely
the full strength of the Schröder-Urysohn Effective Metrization Theorem.
Despite having been used in many different contexts, the full version
of this theorem had not been published before.

Strong computable regularity was introduced by Schröder in \cite{Schroeder1998}
under the name ``computable regularity'', and later renamed by Weihrauch
in \cite{Weihrauch2013}, where other effective notions of regularity
were considered. See Section \ref{sec:Schr=0000F6der's-Metrization-theorem}
for a full definition. 

By establishing that the Effective Metrization Theorem provides an
embedding into the Hilbert cube, we show that the computable second
countability hypothesis is necessary in this theorem. Another characterization
of represented spaces that computably embed into the Hilbert cube
is given in \cite[Theorem 7.1]{Kihara2014}\footnote{This theorem is unfortunately absent from the published version of
\cite{Kihara2014}, namely \cite{Kihara2022}}, in terms of representations computably admitting compact fibers.

The second example that we consider comes from the following classical
statement: 
\begin{equation}
\text{A second countable space is separable. }\label{eq:1-1}
\end{equation}

The natural expected effective version of the above statement uses
the very weak form of effective second countability that we call semi-effective
second countability: 
\begin{prop}
Let $(X,\rho)$ be a semi-effectively second countable represented
space which is overt and has a computable open choice problem. Then
$(X,\rho)$ is effectively separable. 
\end{prop}

However, we prove the following theorem, which implies that the above
proposition is useless:

\begin{theoremalpha}	

A represented space $(X,\rho)$ has a computable open choice problem
if and only if it is computably separable.

\end{theoremalpha}

We then introduce non-total open choice: open choice restricted to
open sets that have a non-empty complement. We then prove:

\begin{theoremalpha}

Having a computable non-total\textbf{ }open\textbf{ }choice does not
imply effective separability. 

\end{theoremalpha}

A correct effective equivalent to statement \eqref{eq:1-1} can be
obtained by using semi-effective second countability and computable
non-total open choice, see Proposition \ref{prop:Non total choice and separability }. 

\subsection{Further notes on the literature }

In Section \ref{sec:Further-notes-on}, we discuss several additional
approaches to computable topology 

Two main points are addressed:
\begin{itemize}
\item The distinction between results that are about spaces up to homeomorphisms
and results about spaces up to computable homeomorphisms.
\item The relationship between the present work and the approach to computable
topology developed by Kalantari and Welch \cite{KW02a,KW03}. In particular,
we present a proposition showing that their approach is not always
compatible with that of Schröder. 
\end{itemize}

\section{\label{sec:Preliminaries}Preliminaries}

The background we require is mostly contained in \cite{Pauly2016}.
See also \cite{Schroeder2021}. To make the paper self contained,
we still introduce the relevant notions. 

\subsection{Represented spaces }

A \emph{representation} of a set $X$ is a partial surjection $\rho:\subseteq\mathbb{N}^{\mathbb{N}}\rightarrow X$.
If $\rho(u)=x$, then $u$ is called a $\rho$\emph{-name} of $x$. 

A point of $X$ is \emph{computable} with respect to $\rho$ if it
admits at least one $\rho$-name which is a computable sequence.

A \emph{multi-function} between sets $X$ and $Y$ is a relation $R\subseteq X\times Y$
which we treat as a function $f$ mapping points of $X$ to subsets
of $Y$ via the formula $f(x)=\{y\in Y\mid(x,y)\in R\}$. We put $\mathrm{dom(}f)=\{x\in X\mid f(x)\neq\emptyset\}$,
and consider that a multi-function $f$ between $X$ and $Y$ is \emph{total}
if $\mathrm{dom(}f)=X$, and that it is \emph{partial} otherwise. 

The fact that $f$ is a total multi-function from $X$ to $Y$ is
denoted via $f:X\rightrightarrows Y$, and we denote a partial multi-function
by $f:\subseteq X\rightrightarrows Y$.
\begin{defn}
Let $(X,\rho)$ and $(Y,\tau)$ be represented spaces, and let $f:\subseteq X\rightrightarrows Y$
be a partial multi-function. A \emph{realizer of $f$ with respect
to $\rho$ and $\tau$, }or $(\rho,\tau)$\emph{-realizer}, is a function
$F:\subseteq\mathbb{N}^{\mathbb{N}}\rightarrow\mathbb{N}^{\mathbb{N}}$
which satisfies the following:
\[
\forall u\in\mathrm{dom}(\rho)\cap\rho^{-1}(\mathrm{dom}(f)),\,\tau(F(u))\in f(\rho(u)).
\]
\end{defn}

In words: $F$ is a realizer of $f$ if it maps the name of any point
of $X$ to a name of one of its images by $f$. 

Let $\rho$ and $\tau$ be representations of a set $X$. Say that
$\rho$\emph{ translates} (resp. \emph{continuously translates}) to
$\tau$ if the identity on $X$ has a computable (resp. continuous)
$(\rho,\tau)$-realizer. 

We denote the fact that $\rho$ translates to $\tau$ by $\rho\le\tau$.
If $\rho\le\tau$ and $\tau\le\rho$, then $\rho$ and $\tau$ are
called\emph{ equivalent, }we define similarly \emph{continuously equivalent
representations}. 
\begin{defn}
[\cite{Kreitz1985,Schroeder2002}]Let $X$ be a topological space.
A representation $\rho:\subseteq\mathbb{N}^{\mathbb{N}}\rightarrow X$
is called \emph{admissible} \emph{for} $X$ if the following hold: 
\begin{enumerate}
\item The representation $\rho$ is continuous, 
\item Every continuous representation $\tau:\subseteq\mathbb{N}^{\mathbb{N}}\rightarrow X$
continuously translates to $\rho$. 
\end{enumerate}
\end{defn}

\begin{thm}
[\cite{Kreitz1985,Schroeder2002}]\label{lem:Admissibility C0 and C0 realizable }Let
$(X,\rho)$ be a represented space and let $(Y,\tau)$ be an admissibly
represented space. Let $f:X\rightarrow Y$ be a function. Then the
following are equivalent: 
\begin{enumerate}
\item The function $f$ is continuous with respect to the final topologies
of $\rho$ and $\tau$.
\item The function $f$ is admits a continuous realizer. 
\end{enumerate}
\end{thm}

\subsection{Representation of continuous functions }

For two represented spaces $(X,\rho)$ and $(Y,\tau)$, there is a
representation, which we denote by $[\rho\rightarrow\tau]$, of all
continuously realizable functions from $(X,\rho)$ to $(Y,\tau)$. 

One way to define this representation is as follows. Continuous functions
on Baire space are exactly the functions that are computable modulo
some oracle. The $[\rho\rightarrow\tau]$-name of a function is thus
the code of a Type 2 machine that computes a realizer, together with
the oracle that this machine requires to function. The fact that currying
and evaluation are computable follows from the Type 2 Universal Turing
Machine Theorem \cite{Weihrauch1985}. For more details, see \cite{Pauly2016}.

\subsection{Sierpiński representation }

The Sierpiński space $\mathbb{S}$ is $\{0,1\}$ with topology generated
by $\{1\}$: $\{1\}$ is open but $\{0\}$ is not. Its standard representation
$c_{\mathbb{S}}:\mathbb{N}^{\mathbb{N}}\rightarrow\mathbb{S}$ is
defined by 
\[
c_{\mathbb{S}}(0^{\omega})=0,
\]
\[
c_{\mathbb{S}}(u)=1\text{ for }u\neq0^{\omega}.
\]

In the present article, we always consider, on any represented space
$(X,\rho)$, that $X$ is equipped with the final topology of $\rho$.
Thus $O$ is open in $X$ if and only if $\rho^{-1}(O)$ is open in
Baire space (or more precisely in $\text{dom}(\rho)$, equipped with
the subset topology of Baire space). 

Because the representation of Sierpiński space is admissible, we get
the following equivalence for any represented space $(X,\rho)$, where
$\mathbf{1}_{O}$ denotes the characteristic function of $O$ : 
\begin{align*}
O\text{ is open in the final topology of }\rho & \iff\mathbf{1}_{O}:X\rightarrow\mathbb{S}\text{ has a continuous }(\rho,c_{\mathbb{S}})\text{-realizer},\\
 & \iff\mathbf{1}_{O}\text{ has a }[\rho\ensuremath{\rightarrow}c_{\mathbb{S}}]\text{-name},
\end{align*}
and thus we can see the representation $[\rho\rightarrow c_{\mathbb{S}}]$
as a representation of the final topology of $\rho$ on $X$. We call
this the \emph{Sierpiński representation associated to} $\rho$. We
denote by $\mathcal{O}(X)$ the represented set thus obtained. 

The computable points of the representation $[\rho\rightarrow c_{\mathbb{S}}]$
are called the\emph{ c.e. open sets. }

The representation $[\rho\rightarrow c_{\mathbb{S}}]$ is also used
to define the represented space $\mathcal{A}_{-}(X)$ of \emph{closed
subsets of $X$ given by negative information: }the name of a closed
set is a $[\rho\rightarrow c_{\mathbb{S}}]$-name of its complement. 

\subsection{Computably admissible representations as $\text{CT}_{0}$ spaces}

Notice that by definition of the Sierpiński representation associated
to a representation $\rho$, the neighborhood map 
\begin{align*}
\mathcal{N}:X & \hookrightarrow\mathcal{O}\mathcal{O}(X)\\
x & \mapsto\{O\mid x\in O\}
\end{align*}
is always computable. 
\begin{thm}
[\cite{Schroeder2002}]A represented space $(X,\rho)$ is admissible
if and only if the neighborhood map $\mathcal{N}$ has a continuously
realizable inverse. 
\end{thm}

Contrary to the original definition, this characterization of admissibility
can be effectivized. This lead Schröder to introducing the following
definition: 
\begin{defn}
[\cite{Schroeder2002}]A representation is \emph{computably admissible}
if the neighborhood map $\mathcal{N}$ has a computable inverse. 
\end{defn}

When $\rho$ is a computably admissible representation of $X$, we
say that the represented space $(X,\rho)$ is \emph{computably Kolmogorov},
or $\text{CT}_{0}$. 

This choice of terminology is based on the following fact: a topological
space is $\text{T}_{0}$ when points are uniquely determined by the
open sets to which they belong. A space is computably Kolmogorov when
the previous equivalence holds computably: the name of a point can
be translated to the name of a program that recognizes the open sets
to which it belongs, and vice versa. 
\begin{thm}
[\cite{Schroeder2002}]A representation $\rho$ of a set $Y$ is
computably admissible if and only if for every represented space $X$
and every function $f:X\rightarrow Y$, we have 
\[
f\text{ is computable}\iff f\text{ is continuous and }f^{-1}:\mathcal{O}(Y)\rightarrow\mathcal{O}(X)\text{ is computable.}
\]
\end{thm}

\subsection{Overt sets }

Let $\textbf{X}=(X,\rho)$ be a represented space.

A subset $A$ of $X$ is called \emph{overt} if the following map
is computable:
\begin{equation}
\mathcal{O}(X)\rightarrow\mathbb{S},\,U\mapsto(U\cap A\ne\emptyset)\label{eq:1-2}
\end{equation}
In other words there is an algorithm that recognizes names of open
sets that meet the set $A$.

A set $A$ is overt if and only if its closure is, and it is only
the closure of $A$ that is uniquely determined by the map given by
(\ref{eq:1-2}).

We denote by $\mathcal{V}(\mathbf{X})$ the \emph{represented space
of closed overt subsets} of $X$: the underlying set is the set of
closed subsets of $X$, equipped with the representation $\psi^{+}$
defined by:
\[
\psi^{+}(p)=A\iff[[\rho\rightarrow c_{\mathbb{S}}]\rightarrow c_{\mathbb{S}}](p)=\{O\in\mathcal{O}(X),\,O\cap A\neq\emptyset\}.
\]
This representation is also often called the\emph{ representation
of closed sets with positive information}. 

The problem Overt Choice was introduced (as \emph{Choice}) in \cite{Brattka1994},
and it has become one of the standard problems in computable analysis
following \cite{Brecht2020}. Hoyrup introduced in \cite{Hoyrup2023}
a variation of it called $\Pi_{2}^{0}$ Overt Choice. The following
problem, Local Overt Choice\footnote{The relevance of this problem, which appears in Figure \ref{fig:Notions-of-bases},
was indicated to us by Arno Pauly.}, is similar to this last problem, but we consider only sets that
are the intersection of a closed set given by overt information with
an open set (given as a characteristic function). 
\begin{defn}
Let $\mathbf{X}$ be a represented space. \emph{Local Overt Choice}
is the following problem: 
\begin{align*}
\mathrm{OVC}_{\mathbf{X}}:\subseteq\mathcal{V}(\mathbf{X})\times\mathcal{O}(\mathbf{X}) & \rightrightarrows X\\
(V,U) & \mapsto V\cap U,
\end{align*}
with $\mathrm{dom(OVC_{\mathbf{X}})}=\{(V,U)\mid V\cap U\neq\emptyset\}$. 
\end{defn}

\subsection{Computably sequential subset}

Let $(X,\rho)$ be a represented space and $Y\subseteq X$ a subset
of it. It is easy to see that the following map is computable:
\begin{align*}
\mathcal{O}(\mathbf{X}) & \rightarrow\mathcal{O}(\mathbf{Y})\\
U & \mapsto U\cap Y
\end{align*}
 By \cite{Schroeder2002}, this map is surjective exactly when the
subset topology on $Y$ is sequential.
\begin{defn}
We say that $Y$ is a \emph{computably sequential} \emph{subset} of
$(X,\rho)$, or that $(Y,\rho_{\vert Y})$ is \emph{computably sequentially
embedded} in $(X,\rho)$, if the projection $\mathcal{O}(\mathbf{X})\rightarrow\mathcal{O}(\mathbf{Y})$
has a computable multivalued inverse. 
\end{defn}

This notion was introduced by Bauer \cite{Bauer2025} in the context
of synthetic topology under the name \emph{intrinsic subset}. 

\section{\label{sec:Notions-of-bases}Notions of bases}

We consider five notions of effective bases/subbases for a represented
space $\textbf{X}=(X,\rho)$. Let $\mathfrak{B}$ be a subset of $\mathcal{O}(X)$,
and $\beta$ a representation of $\mathfrak{B}$. 

\subsection{Semi-effective base}
\begin{defn}
We say that $(\mathfrak{B},\beta)$ is a \emph{semi-effective base
}if the map $\mathfrak{B}\hookrightarrow\mathcal{O}(X)$ is computable
and if $\mathfrak{B}$ is a base of $\mathcal{O}(X)$. 
\end{defn}

Thus the elements of $\mathfrak{B}$ are uniformly open, but the assumption
that $\mathfrak{B}$ is a base is purely a classical one. 

\subsection{Nogina base}
\begin{defn}
We say that $(\mathfrak{B},\beta)$ is a \emph{Nogina base }if it
is a semi-effective base and\emph{ }if furthermore the following multi-function
is computable:
\begin{align*}
N:\subseteq X\times\mathcal{O}(X) & \rightrightarrows\mathfrak{B}\\
(x,O) & \mapsto\{B\in\mathfrak{B},x\in B\subseteq O\}
\end{align*}
Here, $\text{dom}(N)=\{(x,O),x\in O\}\subseteq X\times\mathcal{O}(X)$. 
\end{defn}

This notion of base is closely related to \cite{Nogina1966,Nogina_1969},
this explains our choice of name. It is similar to \cite{Spr98},
however Spreen uses more restrictive conditions based on a \emph{strong
inclusion relation}. This notion is called ``pointwise base'' in
\cite{Bauer2000}. It was also used in \cite{GREGORIADES2016}, in
the computably enumerable version, under the name ``effective countable
base''. 

The following straightforward proposition relates the above definition
to conditions on the base $\mathfrak{B}$ that can be found in \cite{Nogina1966,Nogina_1969}. 
\begin{prop}
If $(\mathfrak{B},\beta)$ is a Nogina base, then the following multifunction
is computable:

\begin{align*}
N':\subseteq X\times\mathfrak{B}\times\mathfrak{B} & \rightrightarrows\mathfrak{B}\\
(x,B_{1},B_{2}) & \mapsto\{B\in\mathfrak{B},x\in B\subseteq B_{1}\cap B_{2}\}
\end{align*}
Here, $\mathrm{dom}(N')=\{(x,B_{1},B_{2}),x\in B_{1}\cap B_{2}\}$. 
\end{prop}

\subsection{Lacombe base}

Consider a semi-effective base $(\mathfrak{B},\beta)$ for $\textbf{X}=(X,\rho)$.
The following fact is well known (see for instance \cite[Corollary 10.2]{Pauly2016},
or \cite{Bauer2012}).
\begin{lem}
\label{lem: Overt unions are open. }The computable injection $i:\mathfrak{B}\hookrightarrow\mathcal{O}(X)$
yields a computable function 
\begin{align*}
j:\begin{cases}
\mathcal{V}(\mathfrak{B}) & \rightarrow\mathcal{O}(X)\\
A & \mapsto\bigcup_{b\in A}i(b)
\end{cases}
\end{align*}
\end{lem}

\begin{proof}
We have the following equivalence, for $x\in X$ and $A\in\mathcal{V}(\mathfrak{B})$:
\begin{align*}
x\in\bigcup_{b\in A}i(b) & \iff\exists b\in A,\,x\in i(b).
\end{align*}
Note that the condition $x\in i(b)$ defines a computable map $X\times\mathfrak{B}\rightarrow\mathbb{S}$,
because $i:\mathfrak{B}\hookrightarrow\mathcal{O}(X)$ is computable.
But, by effective currying, a map of type $X\times\mathfrak{B}\rightarrow\mathbb{S}$
is computable if and only if the corresponding map $X\rightarrow\mathcal{O}(\mathfrak{B})$
is computable. We thus have a computable map $\hat{i}:X\rightarrow\mathcal{O}(\mathfrak{B})$. 

By definition of $\mathcal{V}(\mathfrak{B})$, $\exists:\mathcal{V}(\mathfrak{B})\times\mathcal{O}(\mathfrak{B})\rightarrow\mathbb{S}$
given by $(A,O)\mapsto(\exists x\in A\cap O)$ is computable.

Because $(x,A)\mapsto\exists b\in A,\,x\in i(b)$ is the composition
of $\exists:\mathcal{V}(\mathfrak{B})\times\mathcal{O}(\mathfrak{B})\rightarrow\mathbb{S}$
with $\hat{i}:X\rightarrow\mathcal{O}(\mathfrak{B})$, it is indeed
a computable map on $\mathcal{V}(\mathfrak{B})\times X\rightarrow\mathbb{S}$,
i.e., on $\mathcal{V}(\mathfrak{B})\rightarrow\mathcal{O}(X)$. 
\end{proof}
\begin{defn}
We say that the semi-effective base $(\mathfrak{B},\beta)$ is a \emph{Lacombe
base} if the map $j:\mathcal{V}(\mathfrak{B})\rightarrow\mathcal{O}(X)$
is onto, and if it has a computable multivalued right inverse: a computable
$\psi:\mathcal{O}(X)\rightrightarrows\mathcal{V}(\mathfrak{B})$ such
that $j\circ\psi=\text{id}_{\mathcal{O}(X)}$. 
\end{defn}

In other words, open sets of $X$ can uniformly be written as overt
unions of basic sets. 

The name Lacombe base comes from Lachlan \cite{Lachlan1964} and Moschovakis
\cite{Moschovakis1964} in reference to \cite{Lacombe1957}, where
Lacombe introduced the idea that the computably open sets would be
computable unions of basic open sets. This notion is called ``pointfree
base'' in \cite{Bauer2000}. It was also considered in \cite{Bauer2012}
or in \cite{Brecht2020}. 

Classically, the set of all open sets of a topology is always a base
for this topology. The effective version of this fact is the following
statement.
\begin{prop}
For any represented space, the identity on $\mathcal{O}(X)$ is a
Lacombe base. 
\end{prop}

\begin{proof}
We know by Lemma \ref{lem: Overt unions are open. } that there is
a computable map $\mathcal{V}(\mathcal{O}(X))\rightarrow\mathcal{O}(X)$.
We have to show that it has a computable multivalued inverse. There
is a computable map $\mathcal{O}(X)\rightarrow\mathcal{V}(\mathcal{O}(X))$
given by $U\mapsto\overline{\{U\}}$ \cite[Proposition 7.4(2)]{Pauly2016}.
We show that 
\[
U=\bigcup_{b\in\overline{\{U\}}}b.
\]
This follows directly from the fact that the closure of $\{U\}$ in
the Scott topology is exactly $\{V\in\mathcal{O}(X),\,V\subseteq U\}$.
Indeed, if $V\subseteq U$, then $V\in\overline{\{U\}}$, since Scott
open sets are upper sets. And if $V\not\subseteq U$, then there is
some $x\in V\setminus U$, and $V$ belongs to the Scott open $\mathcal{N}_{x}=\{O\in\mathcal{O}(X),\,x\in O\}$,
while $U$ does not, so $V\notin\overline{\{U\}}$. 
\end{proof}

\subsection{\label{subsec:Representation-subbase}Representation subbase }

Let $X$ be a set, $\mathfrak{B}\subseteq\mathcal{P}(X)$ a set of
subsets of $X$ (which we think of as a subbase for a topology), and
$\beta:\subseteq\mathbb{N}^{\mathbb{N}}\rightarrow\mathfrak{B}$ a
representation of $\mathfrak{B}$. 

Suppose that: 
\begin{enumerate}
\item For every $x$ in $X$, the set $\mathcal{N}_{x}^{\mathfrak{B}}=\{B\in\mathfrak{B}\mid x\in B\}$
is open in the final topology of $\beta$.
\item For every $x$ and $y$ in $X$, $x\neq y$ implies $\mathcal{N}_{x}^{\mathfrak{B}}\neq\mathcal{N}_{y}^{\mathfrak{B}}$.
(Or again: $\mathfrak{B}$ is a subbase for a $\text{T}_{0}$ topology.)
\end{enumerate}
Then we can define a representation of $X$, which we denote by $\beta^{*}$,
by the following formula: 
\[
\forall p\in\mathbb{N}^{\mathbb{N}},\,\beta^{*}(p)=x\iff[\beta\to c_{\mathbb{S}}](p)=\mathcal{N}_{x}^{\mathfrak{B}}.
\]
In other words: the $\beta^{*}$-name of a point $x$ is a name of
the open set $\mathcal{N}_{x}^{\mathfrak{B}}$. Or again: a point
is described its properties. 

The representation $\beta^{*}$ is called the \emph{subbase representation
associated to }$(\mathfrak{B},\beta)$. 

Building on work of Schröder and relying on a Galois connection implicitly
present in \cite{Schroeder2002}, we prove in \cite{Brattka2025}
the following: 
\begin{thm}
\label{thm:1-2-3-theorem-1}A representation of $X$ is admissible
(resp. computably admissible) if and only if it is continuously equivalent
(resp. computably equivalent) to a representation of the form $\beta^{*}$.
\end{thm}

\begin{defn}
We say that $(\mathfrak{B},\beta)$ is a \emph{representation subbase
of $(X,\rho)$ }if $\mathfrak{B}$ is (classically) a subbase of the
final topology of $\rho$, and if furthermore $\rho\equiv\beta^{*}$. 
\end{defn}

In the computably second countable case, when $\mathfrak{B}$ can
be taken to be a totally numbered base, the assumption that it is
a subbase of the final topology of $\rho$ follows automatically from
the equivalence $\rho\equiv\beta^{*}$ (Corollary \ref{cor:Subbasis rep is a subbasis in countable case}).
But in general, the final topology of $\beta^{*}$ can be strictly
finer than the topology generated by the subbase $\mathcal{B}$. This
is precisely the topic of \cite{Brattka2025}, and we refer the interested
reader to this paper for more details. 

As an immediate consequence of Theorem \ref{thm:1-2-3-theorem-1},
we get: 
\begin{prop}
\label{prop: RPZ subbase iff CT0}A represented space has a representation
subbase if and only if it is computably Kolmogorov. 
\end{prop}

\subsection{Enumeration subbase }

The following notion is a generalization of the ``standard representation''
considered by Weihrauch, Kreitz and Grubba in \cite{Kreitz1985,Weihrauch1987,Weihrauch2000,Weihrauch2009ElementaryCT}. 

The idea that a proper generalization of the standard representation
uses a type of formal inclusion relation goes back to Spreen \cite{Spreen2001},
see also \cite{Rauzy2025}. 

The following construction applies only to second-countable spaces,
but the represented base we consider does not have to be countable. 

Consider a represented base $(\mathfrak{B},\beta)$. Consider a relation
$\prec$ on $\text{dom}(\beta)$ which is a \emph{c.e.\ strong inclusion
relation} \cite{Spr98,Spreen2001}, i.e., a relation on $\text{dom}(\beta)$
which satisfies the following two conditions: 
\begin{enumerate}
\item The relation $\prec$ is transitive and semi-decidable: $\prec:\text{dom}(\beta)\times\text{dom}(\beta)\rightarrow\mathbb{S}$
is computable. 
\item For all $p,q$ in $\text{dom}(\beta)$, $p\prec q\implies\beta(p)\subseteq\beta(q)$.
\end{enumerate}
We will make two further assumptions on the relation $\prec$. A \emph{formal
neighborhood base} for a point $x$ is a subset $N_{x}$ of $\text{dom}(\beta)$
that satisfies the two conditions:
\begin{itemize}
\item $\forall b\in N_{x},\,x\in\beta(b),$
\item $\forall b_{1}\in\text{dom}(\beta),\,x\in\beta(b_{1})\implies\exists b_{2}\in N_{x},\,b_{2}\prec b_{1}.$
\end{itemize}
We will assume: 
\begin{enumerate}
\item [(3)]\setcounter{enumi}{3} Every point of $X$ admits a formal neighborhood
base. 
\end{enumerate}
Note that if $\prec$ is reflexive, points will automatically all
admit formal neighborhood bases. Finally, we add the following condition,
which ensures the $\text{T}_{0}$ condition: 
\begin{enumerate}
\item [(4)]\setcounter{enumi}{4} Points are uniquely determined by their
formal neighborhood bases. 
\end{enumerate}
In other words, if $N$ is a formal neighborhood base of $x$ and
$y$, then $x=y$. 

We then define a representation $\beta^{\prec}$ of $X$ by the following:
\[
\beta^{\prec}(\langle u_{0},u_{1},...\rangle)=x\iff\{u_{i},i\in\mathbb{N}\}\text{ is a formal neighborhood base of }x.
\]
(Here $\langle\cdot\rangle:\Pi_{i\in\mathbb{N}}\mathbb{N}^{\mathbb{N}}\rightarrow\mathbb{N}^{\mathbb{N}}$
is a countable tupling function.)

This representation can be thought of as a generalization of the Cauchy
representation of metric spaces: the name of a point is a list of
names of basic open sets that close in on this point. 
\begin{defn}
We say that $(\mathfrak{B},\beta)$ is an \emph{enumeration subbase
for $(X,\rho)$ }if $\rho\equiv\beta^{\prec}$.
\end{defn}

\begin{example}
Consider $\mathbb{R}$ with the Cauchy representation, equipped with
the base that consists of all open intervals. Consider the representation
of this base that comes from the Cauchy representation of $\mathbb{R}$,
and a strong inclusion relation given by $]x,y[\prec]z,t[\iff x>z\,\&\,y<t$.
Then this base is an enumeration subbase. 

Note that we cannot ask, as it is the case in the standard Weihrauch-Kreitz
representation, to have the name of a point $x$ to be a list of all
names of basic open sets that contain it, since a name has to be a
countable sequence. 
\end{example}

\begin{prop}
\label{prop:RPZenumBasis: basic open sets are uniformly open }If
$(\mathfrak{B},\beta)$ is an enumeration subbase for $(X,\rho)$,
then the natural inclusion $\mathfrak{B}\hookrightarrow\mathcal{O}(X)$
is computable. 
\end{prop}

\begin{proof}
We want to show that from the $\beta^{\prec}$-name of a point $x$
and the $\beta$-name $b_{1}$ of a basic set $B$ it is possible
to semi-decide whether $x\in B$. But $x\in B$ if and only if there
appears in the name of $x$ some $\beta$-name $b_{2}$ with $b_{2}\prec b_{1}$.
This is semi-decidable because $\prec$ is. 
\end{proof}
\begin{prop}
\label{prop: Enumeration subbasis =00003D> classical subbasis }If
$(\mathfrak{B},\beta)$ is an enumeration subbase for $(X,\rho)$,
then $\mathfrak{B}$ is (classically) a subbase for $\mathcal{O}(X)$. 
\end{prop}

\begin{proof}
By Proposition \ref{prop:RPZenumBasis: basic open sets are uniformly open },
we know that the elements of $\mathfrak{B}$ are open. Let $x\in X$
and $O\in\mathcal{O}(X)$, with $x\in O$. Then a finite prefix of
the name of $x$ already determines that $x\in O$, since the membership
relation is open in $X\times\mathcal{O}(X)$.

This prefix intersects, via the tupling function $\langle\cdot\rangle$,
finitely many $\beta$-names: $u_{1}$,...,$u_{k}$. Thus $x$ belongs
to the finite intersection $\beta(u_{1})\cap...\cap\beta(u_{n})$,
and we must have $\beta(u_{1})\cap...\cap\beta(u_{n})\subseteq O$,
and $\mathfrak{B}$ is indeed a subbase. 
\end{proof}
A useful property of representations coming from enumeration bases
is that they are open: 
\begin{prop}
\label{prop:The-representation- is open }The representation $\beta^{\prec}$
is open. 
\end{prop}

\begin{proof}
Let $w\in\mathbb{N}^{*}$, we consider the set $\beta^{\prec}(w\mathbb{N}^{\mathbb{N}})$. 

The representation $\beta^{\prec}$ is defined thanks to a tupling
function $\langle\cdot\rangle$. Notice that there exists a tuple
$(w_{1},...,w_{k})$ of elements of $\mathbb{N}^{*}$ such that 
\[
w\mathbb{N}^{\mathbb{N}}=\langle w_{1}\mathbb{N}^{\mathbb{N}},w_{2}\mathbb{N}^{\mathbb{N}},...,w_{k}\mathbb{N}^{\mathbb{N}},\mathbb{N}^{\mathbb{N}},\mathbb{N}^{\mathbb{N}}...\rangle.
\]
Then we have 
\[
\beta^{\prec}(w\mathbb{N}^{\mathbb{N}})=\bigcap_{i\le k}\bigcup_{p\in w_{i}\mathbb{N}^{\mathbb{N}}\cap\text{dom}(\beta)}\beta(p),
\]
which is a finite intersection of unions of open sets, thus open. 
\end{proof}

\section{Relation between the different notions of bases }

Here, we establish the implication relations that relate the different
notions of (sub)bases introduced in Section \ref{sec:Notions-of-bases}.
In Section \ref{sec:Counterexamples-to-separate-bases}, we provide
counterexamples which show that no other relations hold. 

Note that implications between notions of subbases and notions of
bases are obtained by replacing the subbase by the base it induces,
and using the naturally induced representation.
\begin{defn}
If $(\mathfrak{B},\beta)$ is a represented subbase, define the \emph{induced
represented base} $(\cap\mathfrak{B},\cap\beta)$ by:
\begin{itemize}
\item $\cap\mathfrak{B}$ is the set of finite intersections of elements
of $\mathfrak{B}$,
\item The representation $\cap\beta$ is given by:
\[
\text{dom}(\cap\beta)=\{\langle k,u_{0},...,u_{k}\rangle\ensuremath{\in}\mathbb{N}^{\mathbb{N}},k\in\mathbb{N},\,\forall i\le k,\,u_{i}\in\text{dom}(\beta)\};
\]
\[
\forall\langle k,u_{0},...,u_{k}\rangle\in\text{dom}(\cap\beta),\,\cap\beta(\langle u_{0},...,u_{k}\rangle)=\bigcap_{0\le n\le k}\beta(u_{n}).
\]
\end{itemize}
\end{defn}

If $\prec$ was a c.e.\ strong inclusion for $\beta$, we define
the induced strong inclusion on $\text{dom}(\cap\beta)$ as follows:
\[
\langle u_{0},...,u_{k}\rangle\prec\langle v_{0},...,v_{k'}\rangle\iff\forall i\le k',\exists j\le k,\,u_{j}\prec v_{i}.
\]
It is easy to check that this relation remains c.e., and that all
points still admit formal neighborhood bases if this was already the
case. 

The following shows that the operation of replacing a subbase by a
base is harmless from the point of view of the notions of effective
(sub)bases given in Section \ref{sec:Notions-of-bases}. Its proof
is straightforward and left to the reader. 
\begin{prop}
Let $(\mathfrak{B},\beta)$ be a represented subbase for $\mathbf{X}$,
and $(\cap\mathfrak{B},\cap\beta)$ the base it induces. Then $\beta^{*}\equiv(\cap\beta)^{*}$
and $\beta^{\prec}\equiv(\cap\beta)^{\prec}$. 
\end{prop}

We now proceed to prove the part of Theorem \ref{thm:A-All-the-implications between bases}
that concerns implications between notions of represented bases. 

The implications Nogina base $\implies$ Semi-effective base$\implies$
Classical base are straightforward. 

The implication Representation subbase $\implies$ Semi-effective
base is also trivial (modulo taking the induced base), since the definition
of the representation $\beta^{*}$ immediately guarantees that the
map $\mathfrak{B}\rightarrow\mathcal{O}(X)$ will be computable. 
\begin{prop}
An enumeration subbase yields a Nogina base. 
\end{prop}

\begin{proof}
Given a pair $(x,U)$, with $x\in U$, run the program that defines
$U$ with, as input, the name of $x$. This must terminate. Once this
computation has ended, only a finite portion of the name of $x$ was
visited: this finite portion contains the beginning of the names of
finitely many basic open sets whose intersection will serve as a witness
to the Nogina condition. 
\end{proof}
\begin{prop}
A Lacombe base of a computably Kolmogorov space is a representation
subbase. 
\end{prop}

\begin{proof}
Let $\mathfrak{B}$ be the base and denote by $i:\mathfrak{B}\hookrightarrow\mathcal{O}(X)$
the associated computable injection. 

We know that the name of $\mathcal{N}_{x}=\{O\in\mathcal{O}(X),\,x\in O\}\in\mathcal{O}(\mathcal{O}(X))$
can be converted into a $\rho$-name of $x$. 

We have to show that a name of $\mathcal{N}_{x}^{\mathfrak{B}}=\{B\in\mathfrak{B},\,x\in i(B)\}\in\mathcal{O}(\mathfrak{B})$
can also be converted into a $\rho$-name of $x$, this gives a priori
less information, since $i(\mathfrak{B})$ does not contain all open
sets.

But by hypothesis the Sierpiński representation of $\mathcal{O}(X)$
is equivalent to the representation associated to overt unions of
basic sets. 

For an overt $A\subseteq\mathfrak{B}$, $x\in\bigcup_{b\in A}i(b)\iff\exists b\in A,\,x\in i(b)\iff\exists b\in A,\,b\in i^{-1}(\mathcal{N}_{x})\iff\exists b\in A\cap\mathcal{N}_{x}^{\mathfrak{B}}$.
This last condition defines a computable map $\mathcal{V}(\mathfrak{B})\times\mathcal{O}(\mathfrak{B})\rightarrow\mathbb{S}$.
Thus by the smn theorem a name of $\mathcal{N}_{x}^{\mathfrak{B}}\in\mathcal{O}(\mathfrak{B})$
can be translated into a name of the map 
\begin{align*}
\mathcal{V}(\mathfrak{B}) & \rightarrow\mathbb{S}\\
A & \mapsto x\in\bigcup_{b\in A}i(b).
\end{align*}
This is what was to be shown. 
\end{proof}
\begin{prop}
A Lacombe base with computable Local Overt Choice is a Nogina base. 
\end{prop}

\begin{proof}
Let $\mathfrak{B}$ be the base and denote by $i:\mathfrak{B}\hookrightarrow\mathcal{O}(X)$
the associated computable injection. Let $\mathrm{OVC}_{\mathfrak{B}}$
be the local overt choice function associated to $\mathfrak{B}$:
given an overt set and an open set that intersect, it produces a point
in their intersection. 

For $x\in X$, we define $\mathcal{N}_{x}^{\mathfrak{B}}=\{B\in\mathfrak{B}\mid x\in B\}$.
The map $x\mapsto\mathcal{N}_{x}^{\mathfrak{B}}$ is easily seen to
be computable. The Nogina condition is then simply expressed as the
composition of the multifunction $\mathcal{O}(X)\rightrightarrows\mathcal{V}(\mathfrak{B})$
which allows us to express open sets of $X$ as overt unions of basic
set, composed with the following multi-function: 
\begin{align*}
X\times\mathcal{V}(\mathfrak{B}) & \rightrightarrows\mathfrak{B}\\
(x,A) & \mapsto\mathrm{OVC}_{\mathfrak{B}}(A,\mathcal{N}_{x}^{\mathfrak{B}}).
\end{align*}
Both multi-functions are, by hypothesis, computable. 
\end{proof}

\section{\label{sec:Counterexamples-to-separate-bases}Counterexamples that
separate notions of represented bases }

In Section \ref{sec:Counterexamples-to-separate-CE-bases-1}, we provide
counterexamples that separate notions of c.e.\ bases. Thus it is
sufficient here to separate notions of bases that are non-equivalent
in general, but that become equivalent in the case of c.e.\ bases.

\subsection{\label{subsec:A-representation-subbase not Nogina}A representation
subbase that is not a Nogina base }

We now give an example of a representation subbase that does not satisfy
the Nogina condition. In particular, it cannot be an enumeration subbase. 

Consider $X=K\oplus K^{c}=\{2k,k\in K\}\cup\{2k+1,k\notin K\}$. Take
$\mathfrak{B}=\{\{n\},n\in X\}$ with the numbering $\beta:\subseteq\mathbb{N}\rightarrow\mathfrak{B}$
defined by $\beta(n)=\{n\}$, for $n\in X$. 

We consider the subbase representation $\beta^{*}$ associated to
$\beta$: the $\beta^{*}$-name of a point $n$ is a Sierpiński name
of $\{\{n\}\}$. 

Recall that an explicit model of the Sierpiński representation associated
to a representation $\rho$ is as follows \cite{Pauly2016}: a sequence
$(u_{n})_{n\in\mathbb{N}}\in\mathbb{N}^{\mathbb{N}}$ is the name
of an open set $A$ if and only if $u_{0}$ is the code of a Type
2 machine that stops exactly on $\rho$-names of elements of $A$,
and $(u_{n})_{n\ge1}$ is the oracle that this machine requires to
operate. 

For each $i\in\mathbb{N}$, consider the code $t_{i}$ for a Type
2 machine that accepts exactly $2i$ and $2i+1$, and requires no
oracle to function. Then $t_{i}0^{\omega}$ (the concatenation of
$t_{i}$ with the constant zero sequence) is a valid $\beta^{*}$-name,
representing either $2i$ or $2i+1$, depending on whether $i\in K$
or $i\in K^{c}$. 

Suppose that the Nogina condition holds for $\beta^{*}$ and $\beta$,
and let $N$ be a computable realizer of the Nogina condition: if
$n_{x}$ is the name of a point $x$ in $X$ and $n_{A}$ is the name
of some open set $A$ with $x\in A$, then $N(n_{x},n_{A})$ is the
$\beta$-name of a basic set $\{n\}$ with $x\in\{n\}\subseteq A$.
And so $n=x$, and also $N(n_{x},n_{A})=x$. 

Let $n_{X}$ be a computable Sierpiński name of $X$. 

The map 
\begin{align*}
\tilde{N}:\mathbb{N} & \rightarrow\mathbb{N}\\
i & \mapsto N(t_{i}0^{\omega},n_{X})
\end{align*}
 is then computable, and for all $i\in\mathbb{N}$, we have $\tilde{N}(i)=2i$
if $i\in K$, and $\tilde{N}(i)=2i+1$ if $i\in K^{c}$. This is a
contradiction, and the Nogina condition cannot hold for $\beta^{*}$
and $\beta$. 

\subsection{A representation subbase that is not a Lacombe base}

Recall (Theorem \ref{thm:A-All-the-implications between bases}) that
a Lacombe base $\mathfrak{B}$ which admits a computable local overt
choice is automatically a Nogina base. We show that the example given
above (see Section \ref{subsec:A-representation-subbase not Nogina})
does admit a computable local overt choice. Because we have shown
that this base is not a Nogina base, this will imply that it is not
a Lacombe base either. 

The base was $\mathfrak{B}=\{\{n\},n\in X\}$ equipped with the numbering
$\beta:\subseteq\mathbb{N}\rightarrow\mathfrak{B}$ defined by $\beta(n)=\{n\}$,
for $n\in X$. 

What we will show is that the representation of overt subsets of $\mathfrak{B}$
is equivalent to the representation where a set is given by an enumeration
(with a pause symbol) of its elements. It is clear that this representation
allows for a computable local overt choice, and thus this is sufficient
to conclude. 

The name of an overt subset $A$ of $\mathfrak{B}$ encodes a program
which, given an open set $U$ of $\mathfrak{B}$, stops if and only
$U$ intersects $A$. 

Notice that, for each $n\in\mathbb{N}$, the program that accepts
exactly $n$ defines an open subset $U_{n}$ of $\mathfrak{B}$: it
is empty if $n\notin X$, and $\{n\}$ otherwise. 

Thus given an overt set $A$, it is possible to list those $U_{n}$
which intersect $A$, this will precisely give an enumeration of the
elements of $A$.

\subsection{An enumeration subbase that does not induce a Lacombe base }

Here we show that an enumeration subbase does not automatically yield
a Lacombe base. 

In the same setting as above, we will consider a different representation. 

Consider a non-c.e.\ subset $A$ of $\mathbb{N}$. Consider the numbering
$\beta:\subseteq\mathbb{N}\rightarrow\mathfrak{B}$ defined by $\beta(n)=\{n\}$,
for $n\in A$. 

The representation $\beta^{\subseteq}$ that follows from the definition
of an enumeration subbase by taking actual inclusion as a formal inclusion
is equivalent to the natural numbering $\rho$ of $A$ induced by
the identity on $\mathbb{N}$:
\[
\text{dom}(\rho)=A;
\]
\[
\forall n\in\text{dom}(\rho),\,\rho(n)=n.
\]
The $\rho$-c.e.\ open sets are just sets of the form $A\cap E$,
where $E$ is an arbitrary c.e.\ subset of $\mathbb{N}$ (indeed,
$A$ is automatically a computably sequential subset of $\mathbb{N}$).
Yet the computable unions of basic sets are \emph{c.e}.\  \emph{subsets}
of $A$, these form in general a strict subset of the sets of the
form $A\cap E$, for $E$ c.e. For instance take $A\subseteq\mathbb{N}$
to be $K^{c}\oplus K$. Then $A\cap2\mathbb{N}$ is not a c.e.\ subset
of $A$ (because it is not c.e.), yet it is the intersection of a
c.e.\ set with $A$. 

\subsection{An enumeration  subbase which is not a representation subbase }

Putting together the two examples above works: we have a single numbered
base which gives different representations by taking either the enumeration
definition or the representation subbase definition. 

\section{Notions of effective second countability }

\subsection{Main results on computably second countable spaces }

For each of the notions of effective base described in Section \ref{sec:Notions-of-bases},
we get a notion of effective second countability, by replacing the
represented base $\beta:\subseteq\mathbb{N}^{\mathbb{N}}\rightarrow\mathfrak{B}$
by a base equipped with a total numbering: $\beta:\mathbb{N}\rightarrow\mathfrak{B}$.
We will see that notions of bases that were not equivalent become
equivalent.

Note first that the general notion of Lacombe base relies on overt
unions. But on $\mathbb{N}$, the representations of overt subsets
and of open sets coincide with the ``enumeration representation'',
where a set is given by an enumeration that has access to a pause
symbol (so that the empty set can be given by an enumeration). This
simplifies the handling of this representation. 

We now state our theorem on computably second countable spaces. 
\begin{thm}
\label{thm:Strong Effective Second Countability }The following are
equivalent:
\begin{enumerate}
\item $(X,\rho)$ has a totally numbered enumeration subbase.
\item $(X,\rho)$ is a represented space with a representation $\rho$ which
is computably admissible and equivalent to a computably open representation.
\item $(X,\rho)$ has a totally numbered representation subbase.
\item $(X,\rho)$ has a totally numbered Lacombe base and is computably
Kolmogorov.
\end{enumerate}
\end{thm}

Note that the implication $(1)\implies(4)$ can be seen as a Type
2 equivalent of a theorem of Moschovakis \cite{Moschovakis1964} set
in Markovian computable analysis, and which proves that the Nogina
and Lacombe approaches agree on recursive Polish spaces. The hypotheses
required to prove the Type 1 and Type 2 theorems are different: computable
separability is central in \cite{Moschovakis1964} (see the proof
given in \cite{Rauzy2023a}) whereas it plays no role here. 
\begin{proof}
$(1)\implies(2)$ Let $\beta^{\prec}$ be the representation associated
to an enumeration subbase $(B_{n})_{n\in\mathbb{N}}$. The image of
$u_{0}u_{1}...u_{k}\mathbb{N}^{\mathbb{N}}$ by $\beta^{\prec}$ is
the set $B_{u_{0}}\cap...\cap B_{u_{k}}$. A name of this set as an
element of $\mathcal{O}(X)$ can be computed from $u_{0}u_{1}...u_{k}$
by Proposition \ref{prop:RPZenumBasis: basic open sets are uniformly open }.
Thus $\beta^{\prec}$ is computably open. 

We now show that $\beta^{\prec}$ is computably admissible. Notice
that the subbase representation associated to $(B_{n})_{n\in\mathbb{N}}$
translates to $\beta^{\prec}$: for any $x\in X$, given a program
that halts exactly on those $i$ such that $x\in B_{i}$, it is possible
to enumerate the set $\{i\in\mathbb{N},\,x\in B_{i}\}$, this enumeration
yields a name of $x$ for $\beta^{\prec}$. 

Conversely, by Proposition \ref{prop:RPZenumBasis: basic open sets are uniformly open },
the map $\in:X\times\mathbb{N}\rightarrow\mathbb{S},\,(x,n)\mapsto x\in B_{n}$
is computable. This implies that, given a $\beta^{\prec}$-name for
a point $x$, it is possible to compute a name of this point with
respect to the subbase representation associated to $(B_{n})_{n\in\mathbb{N}}$. 

Thus the representation $\beta^{\prec}$ is equivalent to the subbase
representation associated to $(B_{n})_{n\in\mathbb{N}}$. And this
representation yields a $\text{CT}_{0}$ space by Theorem \ref{thm:1-2-3-theorem-1}. 

$(2)\implies(3)$ Suppose that $\rho$ is computably open. The set
$B=\{\rho(w\mathbb{N}^{\mathbb{N}}),w\in\mathbb{N}^{*}\}$, equipped
with its natural numbering obtained via a numbering of $\mathbb{N}^{*}$,
is a totally numbered semi-effective base of $X$. We show that it
is even a representation subbase. 

We must show that a Sierpiński name of the set 
\[
\mathcal{N}_{x}^{B}=\{w\in\mathbb{N}^{*},\,x\in\rho(w\mathbb{N}^{\mathbb{N}})\}\in\mathcal{O}(\mathbb{N}^{*})
\]
can be translated into a $\rho$-name of $x$.

By assumption $\rho$ is computably admissible, and thus any Sierpiński
name of $\mathcal{N}_{x}^{\mathcal{O}(X)}=\{O\in\mathcal{O}(X),\,x\in O\}$
can be translated into a $\rho$-name of $x$. 

It thus suffices to show that a name of $\mathcal{N}_{x}^{B}$ can
be translated to a name of $\mathcal{N}_{x}^{\mathcal{O}(X)}$.

By currying, this is equivalent to showing that the map 
\begin{align*}
\Theta:\subseteq\mathcal{O}(\mathbb{N}^{*})\times\mathcal{O}(X) & \rightarrow\mathbb{S}\\
(\mathcal{N}_{x}^{B},O) & \mapsto x\in O
\end{align*}
is computable (here $\text{dom}(\Theta)=\{U\in\mathcal{O}(\mathbb{N}^{*}),\,\exists x\in X,\,U=\mathcal{N}_{x}^{B}\}\times\mathcal{O}(X)$).
Fix $x\in X$ and $O\in\mathcal{O}(X)$. 

Given a name of the set $\mathcal{N}_{x}^{B}$, it is possible to
computably enumerate all finite prefixes of $\rho$-names of $x$,
denote by $(w_{i})_{i\in\mathbb{N}}$ the obtained sequence. 

The Sierpiński name of $O$ can be seen as a Type 2 machine that accepts
$\rho$-names of points of $O$. If $x$ belongs to $O$, for any
$\rho$-name of $x$, a certain finite prefix of this name is already
accepted by the machine. Therefore 
\[
x\in O\iff\exists i,\,w_{i}\text{ is accepted by the machine encoded in the name of }O.
\]
This is indeed a semi-decidable condition, and $\Theta$ is indeed
computable. 

$(3)\implies(1)$ Take $(B_{n})_{n\in\mathbb{N}}$ a representation
subbase. In this case, the $\rho$-name of a point $x$ is given by
the element $\{n\in\mathbb{N},\,x\in B_{n}\}\in\mathcal{O}(\mathbb{N})$. 

By taking equality of names as a strong inclusion for $(B_{n})_{n\in\mathbb{N}}$,
i.e., $b_{1}\prec b_{2}\iff b_{1}=b_{2}$, we see that the result
follows immediately from the fact that the Sierpiński representation
on $\mathcal{O}(\mathbb{N})$ is equivalent to the following representation:
\begin{align*}
\tau:\mathbb{N}^{\mathbb{N}} & \rightarrow\mathcal{O}(\mathbb{N})\\
p & \mapsto\{p_{n}-1,n\in\mathbb{N}\,\&\,p_{n}>0\}.
\end{align*}

$(1)+(3)\implies(4)$ A representation subbase always yields a computable
Kolmogorov space. 

Replace the subbase by the base it induces, by adding all finite intersections,
and making the intersection map computable. We obtain a numbered base
$(B_{n})_{n\in\mathbb{N}}$.

We show how, given a name of $O\in\mathcal{O}(X)$ with respect to
the Sierpiński representation, we can effectively construct a list
$L$ of names of basic sets such that $\underset{n\in L}{\bigcup}B_{n}=O$. 

We can suppose that the representation $\rho$ of $X$ is an enumeration
representation with respect to $(B_{n})_{n\in\mathbb{N}}$. 

Apply the program $P$ given by the name of $O$ to all finite sequences
$w$, $w\in\mathbb{N}^{*}$. Whenever $w=v_{0}...v_{k}$ is accepted,
add a name of the intersection $B_{v_{0}}\cap...\cap B_{v_{k}}$ to
$L$. 

Remark that an element $w$ is not necessarily the beginning of a
valid name. If it is, and it is accepted by $P$, it must be that
indeed the corresponding intersection $B_{v_{0}}\cap...\cap B_{v_{k}}$
is a subset of $O$. If on the contrary $w$ is not the beginning
of a valid name, it means exactly that $B_{v_{0}}\cap...\cap B_{v_{k}}=\emptyset$,
and thus it will not affect the union that is being constructed. This
shows that $\underset{n\in L}{\bigcup}B_{n}\subseteq O$.

Suppose now that $x\in O$, and take a name $p$ of it. It must be
accepted by $P$, and thus also some finite prefix of it, say $w$.
Thus $w$ defines an intersection $B_{v_{0}}\cap...\cap B_{v_{k}}$
which was added to the list $L$, and thus $x\in\underset{n\in L}{\bigcup}B_{n}$.
This gives the reverse inclusion.

$(4)\implies(3)$ This implication was already true in the general
case of represented bases. 
\end{proof}
By definition, a representation subbase for a represented space $(X,\rho)$
is a classical subbase $\mathcal{B}$ equipped with a representation
$\beta$ such that $\rho\equiv\beta^{*}$. A corollary of the above
proof is that, when considering a totally numbered set $\mathcal{B}$,
the equivalence $\rho\equiv\beta^{*}$ automatically implies that
$\mathcal{B}$ is a subbase of the topology of $X$ (which is the
final topology of $\rho$). 
\begin{cor}
\label{cor:Subbasis rep is a subbasis in countable case}Let $(B_{i})_{i\in\mathbb{N}}$
be a totally numbered set of subsets of $X$, which satisfies the
conditions (1)-(2) of Section \ref{subsec:Representation-subbase},
so that $B$ induces a representation $B^{*}$ of $X$. Then $(B_{i})_{i\in\mathbb{N}}$
is a subbase for the final topology of $B^{*}$.
\end{cor}

\begin{proof}
One can immediately see that the proof of (3)$\implies$(1) of Theorem
\ref{thm:Strong Effective Second Countability } does not use the
assumption that $(B_{i})_{i\in\mathbb{N}}$ is (classically) a subbase
for the topology of $X$. Thus $(B_{i})_{i\in\mathbb{N}}$ is automatically
an enumeration subbase, and Proposition \ref{prop: Enumeration subbasis =00003D> classical subbasis }
shows that an enumeration subbase is automatically a classical subbase. 
\end{proof}
As another corollary, we get: 
\begin{thm}
\label{thm:Intrinsic Embedding}Let $(X,\rho)$ be a computably second
countable represented space, let $Y\subseteq X$ be a subset of $X$,
and equip it with the induced representation $\rho_{\vert Y}$. Then
$(Y,\rho_{\vert Y})$ is also computably second countable, and $Y$
is a computably sequential subset of $(X,\rho)$. 
\end{thm}

\begin{proof}
Let $(B_{n})_{n\in\mathbb{N}}$ be a representation subbase for $(X,\rho)$.
Consider the base $(B_{n}\cap Y)_{n\in\mathbb{N}}$. We claim that
this is a representation subbase of $(Y,\rho_{\vert Y})$. 

For any $y\in Y$, a realizer of the map 
\begin{align*}
\mathcal{N}_{y}^{X}:\mathbb{N} & \rightarrow\mathbb{S}\\
n & \mapsto(y\in B_{n})
\end{align*}
is the same thing as a realizer of 
\begin{align*}
\mathcal{N}_{y}^{Y}:\mathbb{N} & \rightarrow\mathbb{S}\\
n & \mapsto(y\in B_{n}\cap Y).
\end{align*}
Thus a name of $\mathcal{N}_{y}^{Y}$ can be seen as a name of $\mathcal{N}_{y}^{X}$,
which can be turned into a $\rho$-name of $y$ in $X$ (as $(X,\rho)$
is computably second countable), and then seen as a $\rho_{\vert Y}$-name
of $y$. 

And thus $(Y,\rho_{\vert Y})$ is computably second countable. 

By Theorem \ref{thm:Strong Effective Second Countability }, any $O$
in $\mathcal{O}(Y)$ can be written as a countable union 
\[
\bigcup_{n\in A}(B_{n}\cap Y),
\]
with $A\in\mathcal{O}(\mathbb{N})$. This union can also be seen as
the union 
\[
\left(\bigcup_{n\in A}B_{n}\right)\cap Y.
\]
But the union $\bigcup_{n\in A}B_{n}$ is an element of $\mathcal{O}(X)$,
and thus we have uniformly expressed any element of $\mathcal{O}(Y)$
as the intersection of $Y$ with an element of $\mathcal{O}(X)$. 
\end{proof}
The following proposition generalizes the ``computably open representation''
characterization of computable second countability.
\begin{prop}
Let $f:X\rightarrow Y$ be a map between $\mathrm{CT}_{0}$ represented
spaces $X$ and $Y$. If $f$ is computable, computably open and surjective,
and $X$ is computably second-countable, then so is $Y$.
\end{prop}

\begin{proof}
Let $(B_{n})_{n\in\mathbb{N}}$ be a Lacombe base for $X$. We show
that $(f(B_{n}))_{n\in\mathbb{N}}$ is a Lacombe base for $Y$. Because
$f$ is computably open, the map $\mathbb{N}\rightarrow\mathcal{O}(Y),\,n\mapsto f(B_{n})$
is indeed computable. Given an open set $U\subseteq Y$, $f^{-1}(U)$
can be written as a union 
\[
f^{-1}(U)=\bigcup_{i\in I}B_{i}
\]
for a set $I\in\mathcal{O}(\mathbb{N})$ which can be computed from
a name of $U$. Because $f$ is onto, we have $U=f(f^{-1}(U))$, and
so 
\[
U=f(\bigcup_{i\in I}B_{i})=\bigcup_{i\in I}f(B_{i}).
\]
And thus $(f(B_{n}))_{n\in\mathbb{N}}$ is indeed a Lacombe base. 
\end{proof}
The fact that a represented space is computably second countable if
and only if it computably embeds into $\mathcal{O}(\mathbb{N})$ is
well known, it can be found for instance in \cite[Section 2.4.1]{Kihara2019}.
This completes the set of implications about the different notions
of effective second countability that appear in Theorem \ref{thm:A-All-the-implications between bases}. 

Note finally that a representation $\rho$ is computably fiber-overt
if and only if it is computably open \cite{Kihara2014} (a representation
$\rho$ of $X$ is \emph{computably fiber-overt} if the map $X\rightarrow\mathcal{V}(\mathbb{N}^{\mathbb{N}}),\,x\mapsto\rho^{-1}(x)$
is well defined and computable). Computably fiber-overt representations
appear for instance in \cite{Pauly2014,Brattka2018,Brecht2020}. 

\subsection{Additional properties of computably second countable spaces }

In the context of computably second countable spaces, several important
properties can be characterized in terms of properties related to
a base $(B_{n})_{n\in\mathbb{N}}$. 
\begin{prop}
The following are equivalent for a represented space $\mathbf{X}$: 
\begin{enumerate}
\item $\mathbf{X}$ is computably second countable and overt,
\item $\mathbf{X}$ is $\mathrm{CT}_{0}$ and it admits a Lacombe base $(B_{n})_{n\in\mathbb{N}}$
which does not contain the empty set. 
\end{enumerate}
\end{prop}

\begin{proof}
We first prove (1)$\implies$(2). By Theorem \ref{thm:Strong Effective Second Countability },
we can suppose that $\mathbf{X}$ is $\mathrm{CT}_{0}$ and has a
Lacombe base $(B_{n})_{n\in\mathbb{N}}$. But this base could contain
the empty set. However, by overtness, there is a procedure which selects
those basic sets which are not empty. Denote by $\hat{B_{i}}$ the
$i$-th element of $(B_{n})_{n\in\mathbb{N}}$ which is found to be
non-empty by this procedure. We claim that $(\hat{B_{n}})_{n\in\mathbb{N}}$
is also a Lacombe base of $\mathbf{X}$. Indeed, any element of $\mathcal{O}(\mathbf{X})$
can be written as a union of the basic sets $(B_{n})_{n\in\mathbb{N}}$,
and by construction it is immediate to see that such a union can computably
be converted in a union of the basic sets $(\hat{B_{n}})_{n\in\mathbb{N}}$. 

We now prove the converse implication. Suppose that $(B_{n})_{n\in\mathbb{N}}$
is a Lacombe base that does not contain the empty set. In order to
prove that an open set is non-empty, it suffices to write is as the
union of a set of basic sets given by an enumeration. Such an enumeration
uses a special symbol to indicate that nothing is enumerated at certain
stages, and an open set given in this way is non-empty if and only
if a basic set appears in this enumeration, this is indeed semi-decidable.
\end{proof}
Computable points, sets and functions have simple characterizations
in the context of computably second countable spaces: 
\begin{prop}
\label{prop:Characterization of comp points/sets}Let $\mathbf{X}$
be a computably second-countable space with base $(B_{i})_{i\in\mathbb{N}}$.
Then:
\begin{enumerate}
\item $x\in\mathbf{X}$ is computable $\iff\{n\in\mathbb{N}\mid x\in B_{n}\}$
is c.e.,
\item $U\subseteq\mathbf{X}$ is c.e. open $\iff U=\bigcup_{i\in I}B_{i}$
for some c.e. set $I$,
\item $A\subseteq\mathbf{X}$ is overt $\iff\{n\in\mathbb{N}\mid A\cap B_{n}\neq\emptyset\}$
is c.e.,
\item $K\subseteq\mathbf{X}$ is computably compact $\iff\{(n_{1},...,n_{k})\in\mathbb{N}^{*}\mid K\subseteq\bigcup_{i=1}^{k}B_{n_{i}}\}$
is c.e.,
\item For any represented space $\mathbf{Y}$ and function $f:Y\rightarrow X$,
the following are equivalent:
\begin{itemize}
\item $f^{-1}:\mathcal{O}(X)\rightarrow\mathcal{O}(Y),\,U\mapsto f^{-1}(U)$
is computable,
\item $(f^{-1}(B_{n}))_{n\in\mathbb{N}}$ is a computable sequence of c.e.
open sets of $Y$. 
\end{itemize}
\item If $\mathbf{Y}$ is also computably second countable with base $(D_{i})_{i\in\mathbb{N}}$,
the above is also equivalent to:
\begin{itemize}
\item There is a computable function $p:\mathbb{N}^{2}\rightarrow\mathbb{N}$
that satisfies: $\forall n\in\mathbb{N},\,f^{-1}(B_{n})=\bigcup_{i\in\mathbb{N}}D_{p(i,n)}$.
\end{itemize}
\end{enumerate}
\end{prop}

\begin{proof}
All statements are immediate consequences of Theorem \ref{thm:Strong Effective Second Countability }. 
\end{proof}
The interest of Proposition \ref{prop:Characterization of comp points/sets}
lies in the fact that all the right hand side statements have historically
been used as \emph{definitions} (for computable points, computably
open sets, effectively continuous functions, and so on). See for instance
\cite{Lacombe1957,Ceitin1967,Spr98,Iljazovic2018}. Here, however,
we interpret them not as definitions, but as simple characterizations
that hold in the context of computably second-countable spaces. 

If $\mathbf{X}$ is computably second countable, then the represented
space $\mathcal{V}(\mathbf{X})$ of closed overt subsets and the represented
space $\mathcal{K}(\mathbf{X})$ of compact subsets of $\mathbf{X}$
are also computably second countable. 
\begin{thm}
Let $\mathbf{X}$ be a computably second countable represented space
with base $(B_{i})_{i\in\mathbb{N}}$. Then, $\mathcal{V}(\mathbf{X})$
and $\mathcal{K}(\mathbf{X})$ are also computably second countable.
Indeed, the following maps define totally numbered representation
subbases for $\mathcal{V}(\mathbf{X})$ and $\mathcal{K}(\mathbf{X})$
respectively: 
\begin{enumerate}
\item $\diamond B:\mathbb{N}\rightarrow\mathcal{O}(\mathcal{V}(\mathbf{X})),\,n\mapsto\{A\in\mathcal{V}(\mathbf{X})\mid A\cap B_{n}\neq\emptyset\}$. 
\item $\Box B:\mathbb{N}^{*}\rightarrow\mathcal{O}(\mathcal{K}(\mathbf{X})),\,(n_{1},...,n_{k})\mapsto\{K\in\mathcal{K}(\mathbf{X})\mid K\subseteq\bigcup_{i=1}^{k}B_{n_{i}}\}$. 
\end{enumerate}
\end{thm}

The topology on $\mathcal{V}(\mathbf{X})$ is the lower Fell topology
and the topology of $\mathcal{K}(\mathbf{X})$ is the upper Vietoris
topology. 
\begin{proof}
The result follows directly from the following equivalences, together
with Theorem \ref{thm:Strong Effective Second Countability }: 
\begin{itemize}
\item For $A\subseteq X$ closed and $U\subseteq X$ open, with $U=\bigcup_{i\in I}B_{i}$
for some $I\subseteq\mathbb{N}$, we have 
\[
A\cap U\neq\emptyset\iff\exists n\in\mathbb{N},\,A\cap B_{i}\neq\emptyset.
\]
\item For $K\subseteq X$ compact and $U\subseteq X$ open, with $U=\bigcup_{i\in I}B_{i}$
for some $I\subseteq\mathbb{N}$, we have 
\[
K\subseteq U\iff\exists(n_{1},...,n_{k})\in\mathbb{N}^{*},\,A\subseteq\bigcup_{i=1}^{k}B_{n_{i}}.\qedhere
\]
\end{itemize}
\end{proof}
Of course, $\mathbf{X}$ can be computably second countable without
$\mathcal{O}(\mathbf{X})$ being (computably) second countable: this
fails already for $\mathbf{X}=\mathbb{N}^{\mathbb{N}}$, as it is
well known that $\mathcal{O}(\mathbb{N}^{\mathbb{N}})$ is not second
countable. Something else can be said of $\mathcal{O}(\mathbf{X})$
when $\mathbf{X}$ is supposed to be computably second countable:
its representation is computably equivalent to a total representation.
Having a total representation is a form of completeness for represented
spaces (in particular, it coincides with completeness for metric spaces
\cite{Brattka1998} and quasi-metric spaces \cite{Brecht2013}).
\begin{prop}
Let $\mathbf{X}$ be a computably second countable represented space.
Then $\mathcal{O}(\mathbf{X})$ admits a total representation. 
\end{prop}

\begin{proof}
The representation of open sets associated to a Lacombe base $(B_{i})_{i\in\mathbb{N}}$
is total.
\end{proof}
Selivanov studied in \cite{Selivanov2013} total representations of
$\mathcal{O}(X)$ for a second countable space $X$ (see Theorem 8.6
in \cite{Selivanov2013}). But his results do not involve computability. 

\section{\label{sec:Counterexamples-to-separate-CE-bases-1}Counterexamples
that separate notions of effective second countability }
\begin{prop}
The c.e.\ open sets of an admissibly represented space do not have
to generate its topology. 
\end{prop}

\begin{proof}
We consider a total numbering $\nu$ of $\{0,1\}$, given by $\nu(n)=1\iff n\in\text{Tot}$,
where Tot designates $\{n\in\mathbb{N},\text{dom}(\varphi_{n})=\mathbb{N}\}$. 

The numbering $\nu$, seen as a representation, is admissible for
the discrete topology on $\{0,1\}$. However, it is immediate to see
that the only c.e.\ open sets are the empty set and the whole set
itself. 
\end{proof}
An example similar to the above one can be extracted from \cite[Theorem 10]{Hoyrup2016}. 

We now turn to another separation result needed to establish Theorem
\ref{thm:A-All-the-implications between bases}. 
\begin{prop}
\label{prop:The-c.e.-open generate but no semi-eff basis}The c.e.\ open
sets of a $\mathrm{CT}_{0}$ represented space can generate its topology
without it being semi-effectively second countable. 
\end{prop}

We consider an effective version of a well known example of a sequential
but not Fréchet--Urysohn space. This example was studied by Schröder
in \cite{Schroeder2002} as an easy example of an admissibly represented
space which is not second countable. 

We first describe the classical example. Consider the set $X$ given
by:

\[
X=\mathbb{N}^{2}\cup(\{\infty\}\times\mathbb{N})\cup\{(\infty,\infty)\}.
\]

For $m_{0}$, $n_{0}$ in $\mathbb{N}$ and $f:\mathbb{N}\rightarrow\mathbb{N}$,
denote by
\[
D_{n_{0},m_{0}}=\{(n,m_{0}),n\ge n_{0},n\in\mathbb{N}\cup\{\infty\}\};
\]
\[
E_{m_{0},f}=\{(\infty,\infty)\}\cup\bigcup_{m\ge m_{0}}D_{f(m),m}
\]
Define also:
\[
\mathfrak{B}_{1}=\{\{(n,m)\},n,m\in\mathbb{N}\},
\]
\[
\mathfrak{B}_{2}=\{D_{n_{0},m_{0}},n_{0},m_{0}\in\mathbb{N}\},
\]
\[
\mathfrak{B}_{3}=\{E_{m_{0},f},m_{0}\in\mathbb{N},f\in\mathbb{N}^{\mathbb{N}}\}.
\]

Put $\mathfrak{B}=\mathfrak{B}_{1}\cup\mathfrak{B}_{2}\cup\mathfrak{B}_{3}$,
and consider the topology on $X$ generated by $\mathfrak{B}$. 

Note that the set $\mathfrak{B}_{3}$ is not countable. And in $X$,
while $(\infty,\infty)$ is adherent to $\mathbb{N}^{2}$, no sequence
of points of $\mathbb{N}^{2}$ converges to $(\infty,\infty)$, because
this would require that its second component should grow faster that
\emph{any} function $f:\mathbb{N}\rightarrow\mathbb{N}$. 

Here we consider an effective version of this construction: we replace
$\mathfrak{B}_{3}$ by $\mathfrak{B}_{3}^{+}$:
\[
\mathfrak{B}_{3}^{+}=\{E_{m_{0},f},m_{0}\in\mathbb{N},f\in\mathbb{N}^{\mathbb{N}}\text{ is a total computable function}\}
\]
Denote by $\mathfrak{B}^{+}=\mathfrak{B}_{1}\cup\mathfrak{B}_{2}\cup\mathfrak{B}_{3}^{+}$.
We define a numbering $\beta$ of $\mathfrak{B}^{+}$ by 
\[
\beta(\langle1,n,m\rangle)=\{(n,m)\},
\]
\[
\beta(\langle2,n,m\rangle)=D_{n,m},
\]
\[
\beta(\langle3,m,i\rangle)=E_{m,\varphi_{i}}\text{ if \ensuremath{\varphi_{i}} is a total function}.
\]
 To the numbered base $(\mathfrak{B}^{+},\beta)$ we associate the
representation $\rho$ of $X$ which is just the subbase representation.
This guarantees that $(X,\rho)$ is $\mathrm{CT}_{0}$. 

We will now prove that this space is the desired counterexample to
prove Proposition \ref{prop:The-c.e.-open generate but no semi-eff basis}. 

Say that a point $x$ is \emph{effectively adherent }to a set $A$
if the following multi-function is computable: 
\begin{align*}
\Theta_{x,A}:\subseteq\mathcal{O}(X) & \rightrightarrows X\\
O & \mapsto A\cap O,
\end{align*}
with $\text{dom}(\Theta_{x,A})=\{O\in\mathcal{O}(X),\,x\in O\}$.
Notice the following easy result:
\begin{prop}
\label{prop: Eff adherenet limit seq if c.e.  semi-eff basis}If a
represented space $(X,\rho)$ is semi-effectively second countable,
then any point that is effectively adherent to a set is the limit
of a computable sequence of elements of this set. 
\end{prop}

\begin{lem}
\label{lem:The-point-=00005Cinfty not limit comp seq}The point $(\infty,\infty)$,
while effectively adherent to $\mathbb{N}^{2}$, is not the limit
of a computable sequence of points of $\mathbb{N}^{2}$. 
\end{lem}

\begin{proof}
We first show that $(\infty,\infty)$ is effectively adherent to $\mathbb{N}^{2}$.
But this is obvious: any neighborhood of $(\infty,\infty)$ contains
points of $\mathbb{N}^{2}$, one of these can be found by exhaustive
search. 

Now we show that $(\infty,\infty)$ is not the limit of a computable
sequence of points of $\mathbb{N}^{2}$. 

This follows immediately from the fact that if a sequence $((u_{n},v_{n}))_{n\in\mathbb{N}}\subseteq\text{\ensuremath{\mathbb{N}^{2}}}$
converges to $(\infty,\infty)$, then we should have that for every
computable function $f$ and every set $E_{m_{0},f}$, $((u_{n},v_{n}))_{n\in\mathbb{N}}$
eventually belongs to $E_{m_{0},f}$. This implies that eventually
$u_{n}\ge f(v_{n})$ for every computable function $f$. This is not
possible for a computable sequence. 
\end{proof}
Proposition \ref{prop:The-c.e.-open generate but no semi-eff basis}
follows easily. 
\begin{proof}
[Proof of Proposition \ref{prop:The-c.e.-open generate but no semi-eff basis}]
This is simply Proposition \ref{prop: Eff adherenet limit seq if c.e.  semi-eff basis}
together with Lemma \ref{lem:The-point-=00005Cinfty not limit comp seq}.
\end{proof}

\begin{prop}
\label{prop:A-represented-space Lacombe not eff admissible}A represented
space can have a totally numbered Lacombe base without being $\mathrm{CT}_{0}$.
\end{prop}

\begin{proof}
Consider $\{0,1\}$ with the discrete topology and the representation
$\rho:\subseteq\mathbb{N}^{\mathbb{N}}\rightarrow\mathbb{N}^{\mathbb{N}}$
given by 
\[
\rho((u_{n})_{n\in\mathbb{N}})=u_{0}\text{ mod }2,
\]
and $\text{dom}(\rho)=\{u\in\mathbb{N}^{\mathbb{N}},\,u\text{ is not computable}\}$.
This representation is admissible: a translation from the usual representation
$\nu$ of $\{0,1\}$ can be computed by any non-computable oracle. 

We first show that the base given by $B_{0}=\{0\}$ and $B_{1}=\{1\}$
is a Lacombe base for $(\{0,1\},\rho)$.

This follows directly from the fact that the Sierpiński representation
associated to $\rho$ is equivalent to the Sierpiński representation
associated to $\nu$: $[\rho\rightarrow c_{\mathbb{S}}]\equiv[\nu\rightarrow c_{\mathbb{S}}]$. 

To prove this, it suffices to show that the map 
\begin{align*}
\in:\{0,1\}\times\mathcal{O}(\{0,1\}) & \rightarrow\mathbb{S}\\
(x,O) & \mapsto x\in O
\end{align*}
is $(\nu\times[\rho\rightarrow c_{\mathbb{S}}],c_{\mathbb{S}})$-computable
(instead of only $(\rho\times[\rho\rightarrow c_{\mathbb{S}}],c_{\mathbb{S}})$-computable,
as would be expected). 

Consider the $[\rho\rightarrow c_{\mathbb{S}}]$-name of an open set
$O$. It consists of a pair $(u,v)$, where $u\in\mathbb{N}$ is the
code of a Type 2 machine which, when run with access to the oracle
$v\in\mathbb{N}^{\mathbb{N}}$, will accept exactly the $\rho$-names
of points of $O$. Now it is easy to see that $0\in O$ if and only
if some finite sequence $0w\in\mathbb{N}^{*}$ is accepted by this
machine, and $1\in O$ if and only if some finite sequence $1w\in\mathbb{N}^{*}$
is accepted by this machine. These conditions being semi-decidable,
this indeed shows that $\in:\{0,1\}\times\mathcal{O}(\{0,1\})$ is
$(\nu\times[\rho\rightarrow c_{\mathbb{S}}],c_{\mathbb{S}})$-computable. 

Finally, $\rho$ is not computably admissible: if it was, by Theorem
\ref{thm:Strong Effective Second Countability }, the base $(B_{0},B_{1})$
would be a representation subbase, but the representation induced
by the base $(B_{0},B_{1})$ is precisely the natural representation
$\nu$ of $\{0,1\}$, which is not computably equivalent to $\rho$.
\end{proof}
The following proposition contains a represented space which is very
close to being computably second countable -but which still is not. 
\begin{prop}
\label{prop: c.e.  Nogina not Lacombe}A computably admissibly represented
space can be Nogina second countable without being Lacombe second
countable.
\end{prop}

\begin{proof}
We define a representation $\beta$ of a base for the discrete topology
on $\mathbb{N}$ by the following:
\[
\beta(k^{\omega})=\{k\},\,k\in\mathbb{N},
\]
\[
\beta(\ensuremath{\mathbf{1}}_{K})=\text{Tot},
\]
where $\mathbf{1}_{K}$ is the characteristic function of the halting
set. Thus the above is almost the obvious base for the discrete topology
on $\mathbb{N}$, but we artificially add a $\emptyset'$-computable
name to $\text{Tot}$.

Let $\rho$ denote the subbase representation induced by $\beta$.
It is thus computably admissible. 

Notice that $\rho\le\text{id}_{\mathbb{N}}$. The $\rho$-name of
a point $n$ can be thought of as a pair, consisting of $n$ together
with a method to determine whether or not $n$ belongs to $\text{Tot}$
when given access to an oracle for the halting set $K$. 

Let $\tilde{\mathbb{N}}$ be the represented space $(\mathbb{N},\rho)$. 

One immediately checks that the base $\{\{k\},k\in\mathbb{N}\}$ is
a Nogina base for $\tilde{\mathbb{N}}$. Thus $(\mathbb{N},\rho)$
is Nogina second countable. 

We now show that $\tilde{\mathbb{N}}$ does not have a totally numbered
Lacombe base. 

We first consider another base $(B_{i})_{i\in\mathbb{N}}$ for the
discrete topology on $\mathbb{N}$, given by $B_{0}=\text{Tot}$,
$B_{i}=\{i-1\}$ for $i>0$. Consider the associated subbase representation
$\hat{\rho}$. The represented space $(\mathbb{N},\hat{\rho})$ is,
by construction, computably second countable. By Theorem \ref{thm:Strong Effective Second Countability },
the Sierpiński representation $[\hat{\rho}\rightarrow c_{\mathbb{S}}]$
is equivalent to the ``union representation'' associated to $(B_{i})_{i\in\mathbb{N}}$:
open sets can uniformly be written as countable unions of basic sets.
Denote by $\cup B$ this last representation. 

The map $i\mapsto B_{i}$ is a numbering whose image is identical
to the image of $\beta$, it can be seen as a representation which
we denote by $B$. We then have $\beta\le B$, and because the map
$\rho\mapsto[\rho\rightarrow c_{\mathbb{S}}]$ is order reversing
for translation of representations, it follows that $\hat{\rho}\le\rho$
and $[\rho\rightarrow c_{\mathbb{S}}]\le[\hat{\rho}\rightarrow c_{\mathbb{S}}]$. 

Suppose that there exists a Lacombe base $(\mathfrak{B}_{i})_{i\in\mathbb{N}}$
for $\tilde{\mathbb{N}}$, and denote by $\cup\mathfrak{B}$ the representation
of open sets associated to countable unions of elements of $\mathfrak{B}$.
We get: 
\[
\mathfrak{B}\le\cup\mathfrak{B}\equiv[\rho\rightarrow c_{\mathbb{S}}]\le[\hat{\rho}\rightarrow c_{\mathbb{S}}]\equiv\cup B.
\]

By construction of $\rho$, $\text{Tot}$ has an $\emptyset'$-computable
$[\rho\rightarrow c_{\mathbb{S}}]$-name, and thus it also admits
an $\emptyset'$-computable $\cup\mathfrak{B}$-name, which we call
$p_{\text{Tot}}$. Translating this name to $\cup B$, we get a $\emptyset'$-computable
name of $\text{Tot}$ for $\cup B$. But it is clear that a $\cup B$-name
of $\text{Tot}$ either contains explicitly the $B$-name $0$ for
$\text{Tot}$, or it is not $\emptyset'$-computable. Thus the $\emptyset'$-computable
$\cup B$-name of $\text{Tot}$ contains $0$. This implies that a
finite prefix $w_{\text{Tot}}$ of $p_{\text{Tot}}$ is mapped, via
the realizer of the translation $\cup\mathfrak{B}\le\cup B$, to a
finite prefix of a $\cup B$-name that contains $0$. In turn, the
prefix $w_{\text{Tot}}$ can be extended to a computable $\cup\mathfrak{B}$-name
of $\text{Tot}$, for instance the sequence $(w_{\text{Tot}})^{\omega}$
must be a valid $\cup\mathfrak{B}$-name for $\text{Tot}$. 

And thus we get that Tot is actually a c.e.\ open set for $\rho$. 

We conclude by proving that this is impossible. 

For each $i\in\mathbb{N}$, consider a code $t_{i}$ for the program
that, on input $k$, runs $\varphi_{i}(i)$ for $k$ computations
steps, and halts if this computation does not stop during those $k$
steps. Otherwise, it loops indefinitely.

There is a computable function $\Phi$ that transforms $t_{i}$ into
a $\rho$-name of $t_{i}$. Indeed, because, by construction, $t_{i}\in\text{Tot}\iff i\notin K$,
there is an algorithm which, given the $\beta$-name $\ensuremath{\ensuremath{\mathbf{1}}_{K}}$
of $\text{Tot}$, decides whether or not $t_{i}\in\text{Tot}$. 

And thus if $\text{Tot}$ were c.e.\ open for $\rho$, there would
exist a program which, given $i$, would stop if and only if $t_{i}\in\text{Tot}$.
Via the above equivalence this would give a program to semi-decide
the complement of $K$. 
\end{proof}

\section{\label{sec:Schr=0000F6der's-Metrization-theorem}Schröder's Metrization
theorem as a sharp theorem}

Classically, the following equivalences hold for a metric space $(X,d)$:
\begin{enumerate}
\item $(X,d)$ is separable, 
\item $(X,d)$ (topologically) embeds into the Hilbert cube $[0,1]^{\mathbb{N}}$, 
\item $(X,d)$ is second countable. 
\end{enumerate}
Furthermore, the Urysohn Metrization Theorem \cite{Willard1970} shows
that a second countable topological space is metrizable if and only
if it is regular and $\mathrm{T}_{1}$, or equivalently normal and
$\mathrm{T}_{1}$.

The notion of ``effective regularity'' that was found out by Schröder
\cite{Schroeder1998} to be the correct one to establish a metrization
theorem is not the first one one naturally would think of. Different
notions of effective regularity were introduced and compared by Weihrauch
in \cite{Weihrauch2013}. Note however that in \cite{Weihrauch2013}
all spaces are supposed to be computably second countable. 
\begin{defn}
A represented space $(X,\rho)$ is\emph{ computably regular} if the
following multi-function is well defined and computable:
\begin{align*}
R:\subseteq X\times\mathcal{A}_{-}(X) & \rightrightarrows\mathcal{O}(X)^{2}\\
(x,A) & \mapsto\{(U,V),\,x\in U\,\&\,A\subseteq V\,\&\,U\cap V=\emptyset\},
\end{align*}
where $\mathrm{dom}(R)=\{(x,A),\,x\notin A\}$. 

A represented space $(X,\rho)$ is \emph{strongly computably regular}
if the following multi-function is well defined and computable:
\begin{align*}
P:\mathcal{O}(X) & \rightrightarrows\mathcal{O}(X)^{\mathbb{N}}\times\mathcal{A}_{-}(X)^{\mathbb{N}}\\
O & \mapsto\{(U_{n},V_{n})_{n\in\mathbb{N}},\,\forall n\in\mathbb{N},\,U_{n}\subseteq V_{n}\subseteq O,\,O=\bigcup_{n\in\mathbb{N}}U_{n}\}.
\end{align*}

A represented space is \emph{computably normal} if the following multi-function
is well defined and computable:
\begin{align*}
S:\subseteq\mathcal{A}_{-}(X)\times\mathcal{A}_{-}(X) & \rightrightarrows\mathcal{O}(X)\times\mathcal{O}(X)\\
(A,B) & \mapsto\{(U,V),\,A\subseteq U,B\subseteq V,U\cap V=\emptyset\}.
\end{align*}
Here, $\text{dom}(S)=\{(A,B),\,A\cap B=\emptyset\}$. 
\end{defn}

Note the following lemma, which gives sufficient conditions to go
from computable regularity to strong computable regularity:
\begin{lem}
[\cite{Weihrauch2013}]\label{lem:On-overt-and Strong regular}On
overt and computability second countable represented spaces, computable
regularity is equivalent to strong computable regularity. 
\end{lem}

The effective Urysohn lemma states:
\begin{lem}
[Effective Urysohn Lemma, \cite{Schroeder1998}] On a computably
normal space, the following multi-function is computable:
\begin{align*}
R:\subseteq\mathcal{A}_{-}(X)\times\mathcal{A}_{-}(X) & \rightrightarrows\mathcal{C}(X,\mathbb{R})\\
(A,B) & \mapsto\{f,A\subseteq f^{-1}(0),\,B\subseteq f^{-1}(1)\}.
\end{align*}
Here again, $\mathrm{dom}(R)=\{(A,B),A\cap B=\emptyset\}$. 
\end{lem}

Recall that a represented space $(X,\rho)$ \emph{computably embeds}
into a represented space $(Y,\tau)$ if there is a $(\rho,\tau)$-computable
injection $X\hookrightarrow Y$ which admits a $(\tau,\rho)$-computable
partial inverse. 

We now prove the following (which is a slight modification of a result
used in \cite{Amir2023}):
\begin{thm}
[Schröder-Urysohn Effective Metrization]The following are equivalent
for a represented space $(X,\rho)$:
\begin{enumerate}
\item $(X,\rho)$ computably embeds into the Hilbert cube,
\item $(X,\rho)$ computably embeds into some computable metric space,
\item $(X,\rho)$ is computably second countable and strongly computably
regular.
\end{enumerate}
The above imply, but are not equivalent to: 
\begin{enumerate}
\item [(4)]\setcounter{enumi}{4} $(X,\rho)$ is computably second countable
and admits a computable metric that generates its topology. 
\end{enumerate}
\end{thm}

\begin{proof}
$(1)\implies(2)$ is clear. 

$(2)\implies(3)$ Being computably second countable is inherited by
subsets, and so is strong computable regularity.

$(3)\implies(1)$ The first step in the proof given in \cite{Schroeder1998}
is to prove that strong computable regularity and computable second
countability imply computable normality, and thus that the Effective
Urysohn Lemma applies. 

Then, the proof given by Schröder \cite{Schroeder1998} consists in
building a computable double sequence of functions $g_{i,j}:X\to[0,1]$
that separates points, i.e., such that for all $x,y\in X$ with $x\neq y$,
$g_{i,j}(x)\neq g_{i,j}(y)$ for some $(i,j)\in\mathbb{N}^{2}$. 

The functions $g_{i,j}$ are defined as follows. 

Fix $(B_{i})_{i\in\mathbb{N}}$, the countable base that witnesses
computable second countability. By strong computable regularity, there
is a computable double sequence $(U_{i,j},A_{i,j})_{(i,j)\in\mathbb{N}^{2}}\in\mathcal{O}(X)^{\mathbb{N}}\times\mathcal{A}_{-}(X)^{\mathbb{N}}$
such that: 
\[
\forall i\in\mathbb{N},\,B_{i}=\bigcup_{j\in\mathbb{N}}U_{i,j},
\]
\[
\forall(i,j)\in\mathbb{N}^{2},\,U_{i,j}\subseteq A_{i,j}\subseteq B_{i}.
\]

Then, apply the Effective Urysohn’s lemma to $B_{i}^{c}$ and $A_{i,j}$:
this gives a computable function $g_{i,j}$ such that $B_{i}^{c}\subseteq g_{i,j}^{-1}(\{1\})$
and $A_{i,j}\subseteq g_{i,j}^{-1}(\{0\})$. 

In \cite{Schroeder1998} (and in the classical Urysohn theorem), a
metric $d$ is defined as $d(x,y)=\sum_{i,j}2^{-\langle i,j\rangle}|g_{i,j}(x)-g_{i,j}(y)|$. 

Here, as in \cite{Amir2023}, we consider the map $h(x)=(g_{i,j}(x))_{\langle i,j\rangle\in\mathbb{N}}$.
It is by construction a computable map from $(X,\rho)$ to the Hilbert
cube. What we have to show is that it is injective and admits a computable
inverse. 

Let $y$ be a point in $[0,1]^{\mathbb{N}}\cap\text{Im}(h)$. Let
$x=h^{-1}(y)$ be a preimage of $y$. Let $B_{i}$ be some basic set
of $X$.

Notice that by construction, for every $j\in\mathbb{N}$ and $z$
in $B_{i}^{c}$, $g_{i,j}(z)=1$. And thus for every $j\in\mathbb{N}$
and every $z\in X$, $g_{i,j}(z)<1\implies z\in B_{i}$. But, as also
follows from the construction of the functions $g_{i,j}$, it is also
true that for every $z\in B_{i}$, there is some $j\in\mathbb{N}$
such that $g_{i,j}(z)<1$: this holds whenever $z$ belongs to the
open set $U_{i,j}$. 

And thus we get the equivalence: 
\[
x\in B_{i}\iff\exists j\in\mathbb{N},\,g_{i,j}(x)<1.
\]

This equivalence immediately implies what was to be shown:
\begin{itemize}
\item Because $X$ is $\text{T}_{0}$, the above equivalence shows that
the map $h$ is injective. 
\item And because the condition $\exists j\in\mathbb{N},\,g_{i,j}(x)<1$
is semi-decidable in terms of a name of $y=h(x)$, this equivalence
also shows that the map $h$ has a computable inverse. Indeed, it
shows that given the name of an element $y=h(x)$ of the Hilbert cube,
it is possible to compute a name of $\mathcal{N}_{x}^{B}=\{i\in\mathbb{N},\,x\in B_{i}\}\in\mathcal{O}(\mathbb{N})$,
and by computable second countability this name can be translated
into a $\rho$-name of $x$. 
\end{itemize}
Finally, $(1)\implies(4)$ is clear, and the fact that $(4)\centernot\implies(3)$
can be found in \cite[Example 5.4]{Weihrauch2013} (see also \cite{Neu25},
where Weihrauch's example is analyzed in details).
\end{proof}

\section{\label{sec:Open-choice,-non-total}Open choice, non-total open choice,
overtness and separability }

In this section, we study effective versions of the following classical
fact: 
\begin{fact}
\label{fact:A-second-countable-is-Separable }A second countable space
is separable. 
\end{fact}

\begin{proof}
Consider a countable base $(B_{i})$. The set $\{B_{i}\mid B_{i}\neq\emptyset\}$
is also countable. Then apply choice. 
\end{proof}
The effective version of separability is computable separability:
a represented space is \emph{computably separable} if it admits a
dense and computable sequence. 

It is easy to see that in order to obtain an effective version of
the argument above, we will require overtness, to prove that the set
of non-empty subsets of a computably enumerable base is also computably
enumerable. We will also use a form of effective choice axiom, which
we apply only to open sets. We could naively use the following choice
problem: 

\textbf{Open choice: 
\begin{align*}
OC:\mathcal{O}(X)\setminus\{\emptyset\} & \rightrightarrows X\\
O & \mapsto O.
\end{align*}
}

The naive effective version of Fact \ref{fact:A-second-countable-is-Separable }
is the following: 
\begin{prop}
Let $(X,\rho)$ be a semi-effectively second countable represented
space which is overt and has a computable open choice problem. Then
$(X,\rho)$ is effectively separable. 
\end{prop}

The above proposition is obviously true, but it is in fact completely
uninteresting, because of the stronger result, which makes no second
countability assumption: 
\begin{thm}
A represented space $(X,\rho)$ has computable open choice if and
only if it is computably separable. In particular, a space with a
computable open choice is overt. 
\end{thm}

\begin{proof}
The proof relies on the fact that $\{X\}$ is dense in $\mathcal{O}(X)$
for the Scott topology. We first show that this also happens at the
level of names of open sets: if $O$ is an open set of $X$, then
there is a certain name of $O$ which is in the closure of the set
of names of $X$. 

Suppose that $X\neq\emptyset$ (otherwise we have nothing to do).
The subset of $\mathcal{O}(X)$ consisting of $\emptyset$ and of
$X$ is homeomorphic to the Sierpiński space $\mathbb{S}$. It is
easy to see that this is effective: there is a computable embedding
$\mathbb{S}\overset{e}{\hookrightarrow}\mathcal{O}(X)$ mapping $\top$
to $X$ and $\bot$ to $\emptyset$. Denote by $E$ a computable realizer
of $e$. 

It is also well known that the union map $(V,W)\overset{u}{\mapsto}V\cup W$
is computable on $\mathcal{O}(X)$. Let $U$ be the natural computable
realizer of $u$: given as input two open sets $V$ and $W$, each
represented by the code of a Type 2 machine together with the oracle
that machine requires, $U$ produces the code of a machine that halts
if and only if one of these machines halts (and which uses as oracle
the pairing of the two oracles).

For each name $s\in\{0^{n}1^{\omega},n\in\mathbb{N}\}\cup\{0^{\omega}\}$
of an element of the Sierpiński space and each name $p$ of an element
$O$ of $\mathcal{O}(X)$, we consider the name $U(E(s),p)\in\mathbb{N}^{\mathbb{N}}$. 

If $s=0^{n}1^{\omega}$, $U(E(s),p)$ is a name of $X\cup O=X$. If
$s=0^{\omega}$, $U(E(s),p)$ is a name of $\emptyset\cup O=O$. 

The Sierpiński name of an element of $\mathcal{O}(X)$ is an encoded
pair $\langle n,p\rangle$, $n\in\mathbb{N}$ and $p\in\mathbb{N}^{\mathbb{N}}$,
where $n$ is the code for a Type 2 Turing machine and $p$ is the
oracle that this machine will use. 

Let $(w_{n})_{n\in\mathbb{N}}$ be the computable sequence of all
names of the form $\langle n,p\rangle$, for $n\in\mathbb{N}$ and
$p\in\mathbb{N}^{\mathbb{N}}$ an eventually constant sequence. These
elements do not have to be valid names of elements of $\mathcal{O}(X)$,
because we have not guaranteed extensionality: the Turing machine
number $n$ that uses $p$ as oracle can accept some names of a point
and reject others. 

Consider now the computable double sequence $v_{n,t}=U(E(0^{t}1^{\omega}),w_{n})$,
for $n,t\in\mathbb{N}$. Each $v_{n,t}$ is a name of $X$ as an element
of $\mathcal{O}(X)$ -whether or not $w_{n}$ was a valid name of
an element in $\mathcal{O}(X)$. 

Suppose that we have a computable open choice OC for $X$ with computable
realizer $\hat{\text{OC}}$. 

We claim that the computable double sequence $(\rho(\hat{\text{OC}}(v_{n,t})))_{(n,t)\in\mathbb{N}^{2}}$
is dense in $X$. 

First, notice that for each $n,t$, $\hat{\text{OC}}(v_{n,t})$ indeed
belongs to $\text{dom}(\rho)$, since we are applying open choice
to a certain name of $X$, which is non-empty. Now let $O$ be any
non-empty open set of $X$, let $p$ be a name of $O$, and $q=U(E(0^{\omega}),p)$.
Thus $q$ is another name of $O$, and the realizer $\hat{\text{OC}}$
applied to $q$ yields the $\rho$-name of a point $x$ in $O$. By
continuity of the representation $\rho$ of $X$, $\hat{\text{OC}}$
applied to names sufficiently close to $q$ will also yield names
of points in $O$. Such name must appear amongst $(v_{n,t})_{(n,t)\in\mathbb{N}^{2}}$.
\end{proof}
To introduce the correct effectivization of Fact \ref{fact:A-second-countable-is-Separable },
we rely on the following choice problem: 

\textbf{Non-total open choice: 
\begin{align*}
OC^{*}:\mathcal{O}(X)\setminus\{\emptyset,X\} & \rightrightarrows X\\
O & \mapsto O.
\end{align*}
}

We then have: 
\begin{prop}
\label{prop:Non total choice and separability }Let $(X,\rho)$ be
a semi-effectively second countable represented space which is overt
and has a computable non-total open choice problem. Suppose furthermore
that some totally numbered semi-effective base consists only of strict
subsets of $X$. 

Then $(X,\rho)$ is effectively separable. 
\end{prop}

Note that if $X$ can be written as a strict union of two c.e.\ open
sets and is semi-effectively second countable, then there is also
an effective enumeration of a base which never enumerates $X$. Thus
the assumption that some enumerable base avoids $X$ is mild. 

Proposition \ref{prop:Non total choice and separability } is immediate,
but we have to check that it is interesting, by the following: 

\begin{prop}
\label{prop: Non-total OC does not implies Separability}Having computable
non-total\textbf{ }open choice does not imply computable separability,
even on $\mathrm{CT}_{0}$ spaces. 
\end{prop}

\begin{proof}
Consider a set $A\subseteq\mathbb{N}$. Consider the following representation
of $A$: 
\[
\rho(u)=n\iff\text{Im}(u)\cap A=\{n\}.
\]
Thus the $\rho$-name of a point $n$ of $A$ is a list of natural
numbers which intersects $A$ exactly in $n$. 

Note that $\rho$ is computably admissible. Suppose that we have a
name of $\mathcal{N}_{n}^{\mathcal{O}(A)}=\{O\in\mathcal{O}(A),n\in O\}$
as an element of $\mathcal{O}(\mathcal{O}(A))$, i.e., a name that
encodes the characteristic function of $\mathcal{N}_{n}^{\mathcal{O}(A)}$:
on input of the name of an element $U$ of $\mathcal{O}(A)$, it halts
if and only if $n\in U$. We have to recover a $\rho$-name of $n$.
For each $t$ in $\mathbb{N}$, consider the following program $O_{t}$:
on input of a $\rho$-name, it accepts it if and only if it contains
$t$. If $t\in A$, a code for this program is a name of the open
set $\{t\}$ of $\mathcal{O}(A)$. It is not a valid name of an element
of $\mathcal{O}(A)$ when $t\notin A$, because in this case for any
$n\in A$ $O_{t}$ accepts some $\rho$-names of $n$ and rejects
others. When applying the name of $\mathcal{N}_{n}^{\mathcal{O}(A)}$
to the code of $O_{t}$, either $t\in A$, and then $O_{t}$ is accepted
if and only if $t=n$, or $t\notin A$, in which case $O_{t}$ does
not define an element of $\mathcal{O}(A)$ and the behavior of the
realizer of $\mathcal{N}_{n}^{\mathcal{O}(A)}$ on input $O_{t}$
is unspecified: either it never accepts it, or it accepts it. But
in any case, if $O_{t}$ is accepted the realizer of $\mathcal{N}_{n}^{\mathcal{O}(A)}$,
then $t=n$ or $t\notin A$. Thus an enumeration of the set $\{t,O_{t}\text{ is accepted by the realizer of }\mathcal{N}_{n}^{\mathcal{O}(A)}\}$
is a $\rho$-name of $n$: it contains $n$ and possibly some numbers
outside of $A$. 

Furthermore, $(A,\rho)$ has computable non-total\textbf{ }open choice.
Indeed, suppose that some Sierpiński name $P$ of a set $B$ with
$\emptyset\subsetneq B\subsetneq A$ is given. We see $P$ as being
the characteristic function of $B$. We can apply $P$ in parallel
to the following names: $0^{\omega}$, $1^{\omega}$, $2^{\omega}$,...
Notice that if $k\in A$, then $k^{\omega}$ is a valid name of a
point of $A$, and thus it should be accepted at some point if and
only if it belongs to $B$. On the other hand, if $k\notin A$, then
it is an invalid name, but any finite prefix of it could be completed
either into the name of a point of $B$, or into the name of a point
of $A\setminus B$ (which is non-empty by hypothesis). Thus $k$ cannot
be accepted by $P$. And thus the procedure we describe will end up
correctly selecting a point of $B$. 

In the following lemma, Lemma \ref{lem: Non pseudo c.e. set }, we
show that some choice of $A$ guarantees that $(A,\rho)$ is not computably
separable. 
\end{proof}
\begin{lem}
\label{lem: Non pseudo c.e. set }There exits $A\neq\emptyset$ so
that $(A,\rho)$ is not computably separable.
\end{lem}

\begin{proof}
We guarantee that no computable function $g:\mathbb{N}\rightarrow A$
has dense image. Notice that a computable realizer of a computable
function $g:\mathbb{N}\rightarrow A$ is a Type 2 machine that takes
as input a single natural number and outputs a (necessarily computable)
sequence of natural numbers. By the smn-Theorem, we can in fact see
this realizer as a computable function $f:\mathbb{N}\rightarrow\mathbb{N}$
which satisfies that, for each $n$ in $\mathbb{N}$, $W_{f(n)}$
is a $\rho$-name of $g(n)$. (Here $W=\mathrm{dom}(\varphi)$ is
the usual numbering of c.e. subsets of $\mathbb{N}$.) 

Our goal is thus to built $A$ such that there does not exist a total
computable function $f$ such that:
\begin{equation}
\forall n\in\mathbb{N},\,\vert W_{f(n)}\cap A\vert=1;\label{eq:3}
\end{equation}
\[
A\subseteq\bigcup_{n\in\mathbb{N}}W_{f(n)}.
\]
We diagonalize against all computable functions in a non-effective
way.

Suppose that $A\cap\{0,...,N_{i-1}\}$ was already constructed, guaranteeing
that the above conditions do not hold for $\varphi_{k}$, $k=0..i-1$. 

Suppose that $\varphi_{i}$ is total. 

Suppose that for some $n$, $W_{\varphi_{i}(n)}$ contains two points
$x_{1}<x_{2}$ both greater than $N_{i-1}$. Then we add $x_{1}$
and $x_{2}$ to $A$, guaranteeing that \eqref{eq:3} will not be
satisfied for $\varphi_{i}$. We then choose $x_{2}$ to be $N_{i}$:
$A\cap\{0,...,N_{i}\}$ will not change anymore. 

Suppose now that for some $n$, $W_{\varphi_{i}(n)}$ contains a point
already added to $A$ below $N_{i-1}$, and a point $x$ above $N_{i-1}$.
Again we add $x$ to $A$, and \eqref{eq:3} is invalidated. Fix $N_{i}=x$. 

Suppose that for some $n$, $W_{\varphi_{i}(n)}$ contains exactly
one point $x$ above $N_{i-1}$, and possibly some points below $N_{i-1}$,
but that do not belong to $A$. We then add $x+1$ to $A$, and chose
$N_{i}=x+1$. Thus $W_{\varphi_{i}(n)}\cap A$ will be empty. 

Finally, the remaining case is that for every $n$, $W_{\varphi_{i}(n)}$
contains only points below $N_{i-1}$. Then we add $N_{i-1}+1$ to
$A$, and fix $N_{i}=N_{i-1}+1$. 

It is easy to check that all cases are covered, and that the constructed
$A$ is not empty, since each case covered adds at least one point
to $A$. 
\end{proof}

\section{\label{sec:Further-notes-on}Further notes on the literature}

\subsection{Two branches in computable topology }

The study of computable topology in the context of Type 2 computability
can be split in two main branches, which we call \textbf{B1} and \textbf{B2}: 
\begin{description}
\item [{B1}] The study of computable topological spaces up to homeomorphisms,
\item [{B2}] The study of computable topological spaces up to computable
homeomorphisms.
\end{description}
Some surveys and important results in the branch \textbf{B1} are:
\cite{DM23,Sel20a,HMM20,HKS23,KMN25}. Papers illustrating \textbf{B2}
are: \cite{Wei08,Pauly2016,GREGORIADES2016,Brecht2017,Brecht2020a}. 

In terms of represented spaces, these two branches correspond to studying
respectively:
\begin{description}
\item [{B1}] Representations up to continuous translation;
\item [{B2}] Representations up to computable translation.
\end{description}
These two topics are fundamentally different in terms of their goals.
Results that are about spaces up to homeomorphisms are results in
topology, since topology is the domain of mathematics that studies
spaces up to homeomorphism. Thus we could summarize \textbf{B1} as:
\begin{description}
\item [{B1}] Computable topology understood as the application of results
and methods coming from computability in order to study topology.
\end{description}
On the other hand, results about spaces up to computable homeomorphism
are about spaces that have strictly more structure than a topology. 
\begin{description}
\item [{B2}] Computable topology understood as a form of constructive topology,
which provides a refinement of classical results in topology.
\end{description}
Depending on the precise topic studied in computable topology, there
can be more or less overlap between these two branches. 

The present paper's main topic is the study of different notions of
effective second countability. 

On this particular topic, there is in fact absolutely no overlap between
the two branches \textbf{B1} and \textbf{B2} described above. 

Indeed, the following was proved by Hoyrup, Melnikov and Ng:
\begin{thm}
[\cite{Hoyrup2024}]Every second countable space admits a representation
that makes it computably second countable and overt. 
\end{thm}

This theorem shows that: 
\begin{itemize}
\item Up to homeomorphism, every second countable space satisfies the strongest
version of effective second countability. 
\item Weaker notions of effective second countability are relevant only
when spaces are considered up to computable homeomorphism. 
\end{itemize}
And thus the study of the different notions of effective second countability
is a topic in which the two branches \textbf{B1} and \textbf{B2} are
as disjoint as can be. 

But this is not the case of all topics that are studied in computable
topology. 

A nice example, pointed out to us by the referee, is related to the
Birkhoff--Kakutani Theorem, which states that, for second countable
topological groups, Hausdorffness is a sufficient condition for metrizability.

The effective content of this theorem is studied in \cite{KMN25},
where the following is proved:
\begin{thm}
[\cite{KMN25}, Theorem 1.1] For a Polish group $G$ that is either
abelian or locally compact, the following are equivalent:
\begin{enumerate}
\item $G$ has a computable topological presentation,
\item $G$ has a right-c.e.\ Polish presentation.
\end{enumerate}
Furthermore, in (2) the metric can be taken left-invariant. (Or right-invariant.)
\end{thm}

Here, \emph{having a computable topological presentation} is equivalent
to what we call, in the present paper, being computably second countable
and overt, and a \emph{right-c.e.\ presentation of a Polish space}
$X$ is a dense sequence $(u_{n})_{n\in\mathbb{N}}$ such that the
function $\mathbb{N}^{2}\rightarrow\mathbb{R},\,(n,m)\mapsto d(u_{n},u_{m})$
is right computable. 

The subtlety that needs to be pointed out here is the following: 
\begin{itemize}
\item The theorem cited above is a result about spaces up to homeomorphism. 
\item And thus, by looking only at the statement of this theorem, one may
feel that it has no content in terms of the branch \textbf{B2} of
computable topology.
\item However, in order to obtain this result about spaces up to homeomorphism,
an effective version of the Birkhoff--Kakutani Theorem is obtained
in \cite{KMN25}, and this result is valid up to computable homeomorphism. 
\end{itemize}
This illustrates the fact that the distinction between the branches
\textbf{B1} and \textbf{B2} is sometimes thinner than one may expect. 

The effective Birkhoff--Kakutani Theorem, as proved in \cite{KMN25},
recast in the present's paper vocabulary, is the following:
\begin{thm}
[\cite{KMN25}, Theorem 3.2., Effective Birkhoff--Kakutani Theorem]\label{thm:THeorem Melnikov ng}Let
$G$ be a represented group with computable group operations. Suppose
that: 
\begin{itemize}
\item $G$ is computably second countable;
\item Overt;
\item And classically Hausdorff.
\end{itemize}
Then $G$ admits a right-computable left invariant metric which is
compatible with the topological structure of $G$.
\end{thm}

This theorem is very interesting and raises several questions, especially
when compared to the Schröder-Urysohn Metrization Theorem presented
in the present paper. 

We quote three questions which seem to us particularly interesting. 

It was shown by Weihrauch that, by assuming overtness, it is possible
to considerably simplify the hypotheses of the Schröder-Urysohn Metrization
Theorem, because computable regularity can replace strong computable
regularity (see Proposition \ref{lem:On-overt-and Strong regular}).
We thus ask: 
\begin{problem}
Is there a generalization of Theorem \ref{thm:THeorem Melnikov ng}
that does not rely on overtness? 
\end{problem}

Theorem \ref{thm:THeorem Melnikov ng} relies on a non-effective separation
property: the group is assumed Hausdorff, but not necessarily computably
so. As a necessary consequence, the metric that is obtained is only
right-computable. 

On the other hand, in the Schröder-Urysohn Metrization Theorem, the
regularity hypothesis is made effective, and the produced metric is
computable. 
\begin{problem}
Is there an effective Birkhoff--Kakutani Theorem that, by starting
with a computably Hausdorff group, produces a computable metric? 
\begin{problem}
Is there an effective Urysohn Metrization Theorem that, by starting
with a regular but not computably regular space, can characterize
those spaces that admit a right-computable metric? 
\end{problem}

\end{problem}

\subsection{On the work of Kalantari and Welch }

An interesting question, raised by the referee, is to understand how
the present paper relates to the work of Kalantari and Welch on computable
topology, which constitutes a rather important body of work, see for
instance \cite{KW02a,KW03,KW04,KW08,Kalantari2008,KW13,Kalantari2017}. 

This problem has already been addressed several times, in particular
by Kalantari and Welch themselves, and we refer the interested reader
to the last section of \cite{KW03} for a summary of their view of
the problem of integrating their approach into other approaches. 

We want to add the following points: 
\begin{itemize}
\item The approach of Kalantari and Welch to computable topology fits in
the global framework of represented spaces, this is explained in \cite{KW03}. 
\item All the spaces that are considered by Kalantari and Welch are metrizable,
being second countable and regular \cite{KW03}. 
\item A fundamental aspect of the approach of Kalantari and Welch to computable
topology is the use of a certain representation of points associated
to \emph{sharp filters} of basic open sets, see Definition \ref{def:Kalantari-Welch Representation}
below.
\item On the other hand, it seems that Kalantari and Welch have never fixed
a preferred representation of open sets on the spaces that they consider,
and notions such as ``c.e. open sets'' are absent from \cite{KW02a,KW03,KW04,KW08,Kalantari2008,KW13,Kalantari2017}.
\end{itemize}
The representation which is present throughout all the papers of Kalantari
and Welch on computable topology is the following one. 
\begin{defn}
\label{def:Kalantari-Welch Representation} Let $X$ be a second countable
space equipped with a numbered base $\mathfrak{B}=(B_{i})_{i\in\mathbb{N}}$. 
\begin{itemize}
\item A sequence $(B_{u_{n}})_{n\in\mathbb{N}}$ of basic sets defines \emph{a
sharp filter }if:
\begin{enumerate}
\item For all $n$ we have $\overline{B_{u_{n+1}}}\subseteq B_{u_{n}}$;
\item For any two basic sets $C$ and $D$ of $\mathfrak{B}$, with $\overline{C}\subseteq D$,
there is $n$ so that either $B_{u_{n}}\cap C=\emptyset$ or $B_{u_{n}}\subseteq D$. 
\end{enumerate}
\item A sharp filter $(B_{u_{n}})_{n\in\mathbb{N}}$ \emph{converges} to
a point $x$ if 
\[
\{x\}=\bigcap_{n\in\mathbb{N}}B_{u_{n}}.
\]
\item The \emph{Kalantari-Welch representation} associated to $\mathfrak{B}$
is the representation $\delta:\subseteq\mathbb{N}^{\mathbb{N}}\rightarrow X$
given by
\[
\delta(u_{n})=x\iff\text{ (\ensuremath{B_{u_{n}})_{n\in\mathbb{N}}}}\text{ is a sharp filter converging to }x.
\]
\end{itemize}
\end{defn}

The representation of points considered by Kalantari and Welch is
not always compatible with the representations considered in the present
paper. Indeed, we have the following proposition: 
\begin{prop}
There is a computable Polish space in which the Kalantari-Welch representation
of points associated to the basis consisting of rational open balls
gives strictly less information than the Cauchy representation of
points. In fact, the distance function is not right-computable with
respect to this representation. 
\end{prop}

\begin{proof}
We take a certain subset of $\mathbb{R}$. Denote by $K$ the halting
set: $K=\{n\in\mathbb{N},\,\varphi_{n}(n)\downarrow\}$. Consider
the union 
\[
A=\bigcup_{n\in K}[n-\frac{1}{4},n+\frac{1}{4}]\cup\bigcup_{n\notin K}\{n\}.
\]
$A$ is a computably separable subset of $\mathbb{R}$. Indeed, the
set of all natural numbers, together with the set of rationals in
$[n-\frac{1}{4},n+\frac{1}{4}]$, for each $n$ in $K$, constitutes
a computable and dense subset of $A$. Being a computably separable
closed subset of $\mathbb{R}$, it is a computable Polish space. 

The Cauchy representation on $A$ is simply the representation inherited
from the embedding of $A$ in $\mathbb{R}$: a point $x$ is described
by a sequence of rationals $(r_{n})$ with $|x-r_{n}|<2^{-n}$ for
each $n$.

Now we show that the Kalantari-Welch representation associated to
open balls provides strictly less information than the Cauchy representation. 

Suppose that there exists a Type 2 machine $M$ that translates from
the Kalantari-Welch representation to the Cauchy representation. 

For $n\in\mathbb{N}$, we consider the behavior of $M$ if its input
is the constant sequence of open balls $B(n,1/2)$, $B(n,1/2)$, $B(n,1/2)$,
$B(n,1/2)$... 
\begin{itemize}
\item If $n\notin K$, then this is a valid name with respect to the Kalantari-Welch
representation, since $B(n,1/2)=\{n\}$. (Condition (2) of Definition
\ref{def:Kalantari-Welch Representation} is trivially satisfied for
a singleton.) And so $M$ should start writing a Cauchy name of $n$,
in particular it should, at some point, provide an approximation of
$n$ within $1/8$. 
\item If $n\in K$, since we do have $\overline{B(n,1/2)}\subseteq B(n,1/2)$,
any finite sequence $B(n,1/2)$, $B(n,1/2)$, ..., $B(n,1/2)$ constitutes
the beginning of a valid name for \emph{any} point in $[n-\frac{1}{4},n+\frac{1}{4}]$.
And in this case the machine $M$ can never provide an approximation
of $n$ within $1/8$, since this would necessarily exclude numbers
that could later turn out to correspond to the named number.
\end{itemize}
And thus we see that the machine $M$ should produce an approximation
of its input within $1/8$ exactly for those $n$ that do not belong
to $K$. This is a contradiction because the complement of $K$ is
not c.e. Thus the Kalantari-Welch representation does not translate
to the Cauchy representation.

One easily notices that the example above also implies that the distance
fu nction is not right-computable with respect to the Kalantari-Welch
representation. 
\end{proof}
\bibliographystyle{alpha}
\bibliography{../Biblio/TTEMOSCHO}

\end{document}